\documentclass{article}

\usepackage{arxiv}
\usepackage{natbib}

\usepackage[utf8]{inputenc} 
\usepackage[T1]{fontenc}    
\usepackage{hyperref}       
\usepackage{url}            
\usepackage{booktabs}       
\usepackage{amsfonts}       
\usepackage{nicefrac}       
\usepackage{microtype}      
\usepackage{lipsum}		
\usepackage{graphicx}
\usepackage{natbib}
\usepackage{doi}

\usepackage{amsmath,amsbsy,amsgen,amscd,amssymb,bm,xspace}
\usepackage[capitalise]{cleveref}
\usepackage{amsfonts}

\usepackage[noend]{algpseudocode}
\usepackage{algorithm,algorithmicx}
\algrenewcommand\alglinenumber[1]{\sf\tiny\color{medblue}{#1}\quad}
\algrenewcommand\algorithmicrequire{\textbf{input:}}
\algrenewcommand\algorithmicensure{\textbf{output:}}
\algdef{SE}[DOWHILE]{Do}{doWhile}{\algorithmicdo}[1]{\algorithmicwhile\ #1}%
\usepackage{multirow}
\usepackage{booktabs} 
\usepackage{rotating}
\usepackage{subcaption}
\newcommand{\mnote}[1]{}
\newcommand{\snote}[1]{}
\newcommand{\znote}[1]{}
\newcommand{\bnote}[1]{}

\renewcommand{\mnote}[1]{\textcolor{blue}{\textbf{[MU: #1]}}}
\renewcommand{\snote}[1]{\textcolor{red}{\textbf{[SZ: #1]}}}
\renewcommand{\znote}[1]{\textcolor{purple}{\textbf{[ZF: #1]}}}
\renewcommand{\bnote}[1]{\textcolor{cyan}{\textbf{[Bartolomeo: #1]}}}

\usepackage{mathtools}
\usepackage{array}
\usepackage{ragged2e} 
\usepackage{titletoc}
\usepackage{pifont}
\newcommand{\xmark}{\ding{55}}%
\newcommand{\bigO}{\mathcal O}
\newcommand{\redX}{\textcolor{Maroon}{\xmark}}

\newcolumntype{L}[1]{>{\raggedright\let\newline\\\arraybackslash\hspace{0pt}}m{#1}}
\newcolumntype{C}[1]{>{\centering\let\newline\\\arraybackslash\hspace{0pt}}m{#1}}
\newcolumntype{R}[1]{>{\raggedleft\let\newline\\\arraybackslash\hspace{0pt}}m{#1}}

\usepackage{makecell}

\usepackage[dvipsnames]{xcolor}
\graphicspath{{figures/}}
\newcommand{\reals}{\mathbb{R}}

\newcommand{\sym}{\mathbb{S}} 

\usepackage{amsthm}
\usepackage{mathtools,bm}

\newif\ifpreprint
\preprinttrue

\newcommand{\argmin}{\mathop{\rm argmin}}

\newcommand{\lamM}{\lambda_{\textup{min}}(MM^T)}
\newcommand{\Nm}{N^{-}}
\newcommand{\Err}{\textrm{Err}^k_{\nabla}}

\newcommand{\epsTot}{\varepsilon_w^k}
\newcommand{\epsjTot}{\varepsilon^j_{w}}
\newcommand{\theMax}{\theta_{\textup{max}}}
\newcommand{\theMin}{\theta_{\textup{min}}}
\newcommand{\psiMax}{\psi_{\textup{max}}}

\newcommand{\prox}{\textup{prox}}
\newcommand{\epsappx}{\overset{\varepsilon}\approx}

\newcommand{\xepsappx}{\overset{\varepsilon_x^k}\approx}
\newcommand{\zfepsappx}{\overset{\varepsilon_z^k}\approxeq}

\newtheorem{theorem}{Theorem}
\newtheorem{proposition}{Proposition}%
\newtheorem{definition}{Definition}%
\newtheorem{lemma}{Lemma}
\newtheorem{corollary}{Corollary}
\newtheorem{assumption}{Assumption}
\newtheorem{remark}{Remark}

\title{On the (linear) convergence of Generalized Newton Inexact ADMM}

\author{{Zachary Frangella}\\
	Stanford University\\
	\texttt{zfran@stanford.edu} \\
	\And
	{Theo Diamandis} \\
	Massachusetts Institute of Technology\\
	\texttt{tdiamand@mit.edu} \\
        \And
        {Bartolomeo Stellato} \\
      Princeton University \\
      \texttt{bstellato@princeton.edu} \\
      \And 
      {Madeleine Udell} \\
      Stanford University \\ 
      \texttt{udell@stanford.edu}
}

\renewcommand{\shorttitle}{Generalized Newton Inexact ADMM}

\hypersetup{
pdftitle={Generalized Newton Inexact ADMM},
}

\begin{document}
\maketitle

\begin{abstract}
	This paper presents GeNI-ADMM, a framework for large-scale composite convex optimization
that facilitates theoretical analysis of both existing and new approximate ADMM schemes.
GeNI-ADMM encompasses any ADMM algorithm that solves a first- or second-order approximation
to the ADMM subproblem inexactly.
GeNI-ADMM exhibits the usual $\mathcal{O} (1/t)$-convergence rate under standard hypotheses
and converges linearly under additional hypotheses such as strong convexity.
Further, the GeNI-ADMM framework provides explicit convergence rates
for ADMM variants accelerated with randomized linear algebra,
such as NysADMM and sketch-and-solve ADMM,
resolving an important open question on the convergence of these methods.
This analysis quantifies the benefit of improved approximations
and can aid in the design of new ADMM variants with faster convergence.
\end{abstract}


\section{Introduction}
\label{section:intro}
The Alternating Direction Method of Multipliers (ADMM) is one of the most popular methods for solving composite optimization problems, as it provides a general template for a wide swath of problems and converges to an acceptable solution within a moderate number of iterations~\citep{MAL-016}.
Indeed,~\cite{MAL-016} implicitly promulgates the vision that ADMM provides a unified solver for various convex machine learning problems. 
Unfortunately, for the large-scale problem instances routinely encountered in the era of Big Data, ADMM scales poorly and cannot provide a unified machine learning solver for problems of all scales.
The scaling issue arises as ADMM requires solving two subproblems at each iteration, whose cost can increase superlinearly with the problem size.
As a concrete example, in the case of $\ell_1$-logistic regression with an $n\times d$ data matrix, ADMM requires solving an $\ell_2$-regularized logistic regression problem at each iteration \citep{MAL-016}. 
With a fast-gradient method, the total complexity of solving the subproblem is $\tilde {\mathcal{O}}(nd\sqrt{\kappa})$ \citep{bubeck2015convex}, where $\kappa$ is the condition number of the problem.
When $n$ and $d$ are in the tens of thousands or larger---a moderate problem size by contemporary machine learning standards---and $\kappa$ is large, such a high per-iteration cost becomes unacceptable.
Worse, ill-conditioning is ubiquitous in machine learning problems; often $\kappa = \Omega(n)$, in which case the cost of the subproblem solve becomes superlinear in the problem size.  

Randomized numerical linear algebra (RandNLA) offers promising tools to scale 
ADMM to larger problem sizes.
Recently \citet{pmlr-v162-zhao22a} proposed the algorithm NysADMM,
which uses a randomized fast linear system solver to scale ADMM up to problems with tens of thousands of samples and hundreds of thousands of features. 
The results in \citet{pmlr-v162-zhao22a} show that ADMM combined with the RandNLA primitive runs \emph{3 to 15$\times$ faster} than state-of-the-art solvers on 
machine learning problems from LASSO to SVM to logistic regression.
Unfortunately, the convergence of randomized or approximate ADMM solvers like NysADMM
is not well understood.
NysADMM approximates the $x$-subproblem using a linearization 
based on a second-order Taylor expansion, 
which transforms the $x$-subproblem into a Newton-step, {\it i.e.}, a linear system solve.
It then solves this system approximately (and quickly) using a randomized linear system solver.
The convergence of this scheme, which combines linearization and inexactness, 
is not covered by prior theory for approximate ADMM;
prior theory covers either linearization~\citep{ouyang2015accelerated} or inexact solves~\citep{eckstein1992douglas} but not both.

In this work, we bridge the gap between theory and practice to explain the excellent 
performance of a large class of approximate linearized ADMM schemes, including NysADMM~\citep{pmlr-v162-zhao22a} and many methods not previously proposed. 
We introduce a framework called Generalized Newton Inexact ADMM, which we refer to as GeNI-ADMM (pronounced \emph{genie-ADMM}). 
GeNI-ADMM includes NysADMM and many other approximate ADMM schemes as special cases. 
The name is inspired by viewing the linearized $x$-subproblem in GeNI-ADMM 
as a generalized Newton-step.
GeNI-ADMM allows for inexactness in both the $x$-subproblem and the $z$-subproblem. 
We show GeNI-ADMM exhibits the usual $\mathcal O(1/t)$-convergence rate under standard assumptions, 
with linear convergence under additional assumptions.
Our analysis also clarifies the value of using curvature in the generalized Newton step:
approximate ADMM schemes that take advantage of curvature
converge faster than those that do not at a rate that depends on the conditioning of the subproblems.
As the GeNI-ADMM framework covers any approximate ADMM scheme that replaces 
the $x$-subproblem by a linear system solve, 
our convergence theory covers any ADMM scheme that uses fast linear system solvers.
Given the recent flurry of activity on fast linear system solvers within the (randomized) 
numerical linear algebra community \citep{lacotte2020effective,meier2022randomized,frangella2023randomized}, our results will help realize these benefits  
for optimization problems as well.
To demonstrate the power of the GeNI-ADMM framework, we establish convergence of NysADMM and another RandNLA-inspired scheme, sketch-and-solve ADMM, whose convergence was left as an open problem in~\cite{buluc2021randomized}.

\subsection{Contributions}
Our contributions may be summarized concisely as follows:
\begin{enumerate}
    \item We provide a general ADMM framework GeNI-ADMM, that encompasses prior approximate ADMM schemes as well as new ones. 
    It can take advantage of second-order information and allows for inexact subproblem solves.
    \item Our analysis shows the benefits of schemes that employ preconditioning and variable metrics over methods that do not.  
    \item We show that despite all the approximations it makes, GeNI-ADMM converges ergodically at the usual ergodic $\mathcal O(1/t)$ rate.
    \item We apply our framework to show some RandNLA-based approximate ADMM schemes converge at the same rate as vanilla ADMM, answering some open questions regarding their convergence \citep{buluc2021randomized,pmlr-v162-zhao22a}.
    In the case of sketch-and-solve ADMM, we show modifications to the naive scheme are required to ensure convergence.
\end{enumerate}

\subsection{Roadmap} In \cref{section:problemstatement}, 
we formally state the optimization problem that we focus on in this paper and briefly introduce ADMM. 
In \cref{section:GeneralizedNewtADMM} we introduce the Generalized Newton Inexact ADMM framework and review ADMM and its variants. 
\Cref{section:TechPrelimsAssump} gives various technical backgrounds and assumptions needed for our analysis. 
\Cref{section:SubLinConvergence} establishes that GeNI-ADMM converges at an $\mathcal O(1/t)$ rate in the convex setting. 
In \Cref{section:applications}, we apply our theory to establish convergence rates for two methods that naturally fit into our framework, 
and we illustrate these results numerically in \Cref{section:experiments}.

\subsection{Notation and preliminaries}
We call a matrix psd if it is positive semidefinite. We denote the convex cone of $n\times n$ real symmetric psd matrices by $\sym^+_n$. 
We denote the Loewner ordering on $\sym^+_n$ by $\preceq$, that is $A\preceq B$ if and only if $B-A$ is psd.
Given a matrix $H$, we denote its spectral norm by $\|H\|$.
If $f$ is a smooth function we denote its smoothness constant by $L_f$.
We say a positive sequence $\{\varepsilon^{k}\}_{k\ge 1}$ is summable if $\sum_{k=1}^{\infty}\varepsilon^{k}<\infty$.

\section{Problem statement and ADMM}
\label{section:problemstatement}
Let $\mathcal{X}, \mathcal{Z}, \mathcal H$ be finite-dimensional inner-product spaces with inner-product $\langle \cdot,\cdot \rangle$ and norm $\|\cdot\|$.
We wish to solve the convex constrained optimization problem
\begin{equation}
\label{eq:ADMMProb}
\begin{array}{ll}
\mbox{minimize} & f(x)+g(z) \\
\mbox{subject to} & Mx+Nz = 0 \\
& x\in X,~ z\in Z,
\end{array}
\end{equation}
with variables $x$ and $z$, 
where $X\subset \mathcal{X}$ and $Z\subset \mathcal{Z}$ are closed convex sets, 
$f$ is a smooth convex function, 
$g$ is a convex proper lower-semicontinuous (lsc) function, 
and $M: \mathcal{X} \mapsto \mathcal{H}$ and $N:\mathcal{Z}\mapsto\mathcal{H}$ are bounded linear operators.
\begin{remark}
    Often, the constraint is presented as $Mx+Nz = c$ for some non-zero vector $c$ however, by increasing the dimension of $\mathcal Z$ by $1$ and replacing $N$ by $[N~c]$, we can make $c = 0$. Thus, our setting of $c = 0$ is without loss of generality. 
\end{remark}
We can write problem~\eqref{eq:ADMMProb} as the saddle point problem 
\begin{equation}
\label{eq:saddlepoint}
\underset{{x \in X, z \in Z}} {\mbox{minimize}}\; \underset{y \in Y}{\mbox{maximize}} \; f(x) + g(z) +\langle y, Mx+Nz\rangle, 
\end{equation}
where $Y \subset \mathcal Z$ is a closed convex set. 
The saddle-point formulation will play an important role in our analysis. 
Perform the change of variables $u = y/\rho$ and define the Lagrangian
\[
L_\rho(x,z,u) \coloneqq  f(x) + g(z) +\langle \rho u, Mx+Nz \rangle.
\] 
Then \eqref{eq:saddlepoint} may be written concisely as
\begin{equation}
    \underset{x \in X, z \in Z} {\mbox{minimize}}\; \underset{{u \in U}}{\mbox{maximize}} \; L_\rho(x,z,u).
\end{equation}
The ADMM algorithm (\cref{alg-admm}) is a popular method for solving~\eqref{eq:ADMMProb}. 
Our presentation uses the scaled form of ADMM \citep{MAL-016,ryu2022large}, 
which uses the change of variables $u = y / \rho$, 
this simplifies the algorithms and analysis, 
so we maintain this convention throughout the paper. 
\begin{algorithm}[H]
	\centering
	\caption{ADMM\label{alg-admm}}
	\begin{algorithmic}
		\Require{convex proper lsc functions $f$ and $g$, constraint matrix $M$, stepsize $\rho$}
		\Repeat
		\State{$x^{k+1} = \underset{x \in X}{\argmin} \{f(x) + \frac{\rho}{2}\|Mx+Nz^k+ u^k\|^2\}$}
		\State{$z^{k+1} = \underset{z \in Z}{\argmin} \{g(z) + \frac{\rho}{2}\|Mx^{k+1} +Nz +u^k\|^2\}$}
		\State{$u^{k+1} = u^k + Mx^{k+1} +Nz^{k+1} $}
		\Until{convergence}
		\Ensure{solution $(x^\star, z^\star)$ of problem \eqref{eq:ADMMProb}}
	\end{algorithmic}
\end{algorithm} 
\section{Generalized Newton Inexact ADMM}
\label{section:GeneralizedNewtADMM}

As shown in \cref{alg-admm}, at each iteration of ADMM, two subproblems are solved sequentially to update variables $x$ and $z$. 
ADMM is often the method of choice when the $z$-subproblem has a closed-form solution. 
For example, if $g(x) = \|x\|_1$, the $z$-subproblem is soft thresholding, 
and if $g(x) = 1_{\mathcal S}(x)$ is the indicator function of a convex set $\mathcal{S}$, 
the $z$-subproblem is projection onto $\mathcal{S}$
\cite[\S6]{parikh2014proximal}.
However, it may be expensive to compute even with a closed-form solution.
For example, when $g(x) = 1_{\mathcal S}(x)$ and $\mathcal S$ is the psd cone, the $z$-subproblem requires an eigendecomposition to compute the projection, which is prohibitively expensive for large problems \citep{rontsis2022efficient}.

Let us consider the $x$-subproblem
\begin{equation}
    \label{eq:xsubproblem}
    x^{k+1} = \underset{x \in X}{\argmin} \left\{f(x) + \frac{\rho}{2}\|Mx+Nz^k + u^k\|^2\right\}.
\end{equation}
In contrast to the $z$-subproblem, there is usually no closed-form solution for the $x$-subproblem. 
Instead, an iterative scheme is often used to solve it inaccurately,
especially for large-scale applications.
This solve can be very expensive when the problem is large. 
To reduce computational effort,
many authors have suggested replacing this problem with 
a simplified subproblem that is easier to solve.
We highlight several strategies to do so below.

\paragraph{Augmented Lagrangian linearization.}
One strategy is to linearize the augmented Lagrangian 
term $\frac{\rho}{2}\|Mx+Nz^k + u^k\|^2$ in the ADMM subproblem
and replace it by the quadratic penalty $\frac{1}{2}\|x - x^k\|^2_P$ for some 
(carefully chosen) positive definite matrix $P$.
More formally, the strategy adds a quadratic term to form a new subproblem 
\begin{equation*}
    x^{k+1} = \underset{x \in X}{\argmin} \left\{f(x) + \frac{\rho}{2}\|Mx+Nz^k + u^k\|^2 + \frac{1}{2}\|x - x^k\|^2_P\right\},
\end{equation*}
which is substantially easier to solve for an appropriate choice of $P$.
One canonical choice is $P = \eta I - \rho M^TM$, where $\eta > 0$ is a constant.
For this choice, the quadratic terms involving $M$ cancel,
and we may omit constants with respect to $x$,
resulting in the subproblem 
\begin{equation*}
    x^{k+1} = \underset{x \in X}{\argmin} \left\{f(x) 
    + \rho \langle Mx, Mx^k - z^k + u^k \rangle
    + \frac{\eta}{2}\|x - x^k\|^2 \right\}.
\end{equation*}
Here, an isotropic quadratic penalty replaces the augmented Lagrangian term in \eqref{eq:xsubproblem}.
This strategy allows the subproblem solve to be replaced by a proximal operator with a (possibly) closed-form solution. 
This strategy has been variously called \emph{preconditioned ADMM}, proximal ADMM, and (confusingly) linearized ADMM \citep{deng2016global,he20121,ouyang2015accelerated}.

\paragraph{Function approximation.}
The second strategy to simplify the $x$-subproblem, is to approximate the function $f$ by a first- or second-order approximation \citep{ouyang2015accelerated,ryu2022large,pmlr-v162-zhao22a}, forming the new subproblem
\begin{equation}
\label{eq:functionapprox}
    x^{k+1} = \underset{x \in X}{\argmin} \left\{f(x^k) 
    + \langle \nabla f(x^k), x - x^k\rangle 
    + \frac 1 2 \| x - x^k \|_H^2 
    + \frac{\rho}{2}\|Mx+Nz^k + u^k\|^2\right\}
\end{equation}
where $H$ is the Hessian of $f$ at $x^k$.
The resulting subproblem is quadratic and may be solved 
by solving a linear system or (for $M=I$) performing a linear update, 
as detailed below in \cref{subsection:recoveredadmm}. 

\paragraph{Inexact subproblem solve.}
The third strategy is to solve the ADMM subproblems inexactly
to achieve some target accuracy, either in absolute error or relative error. 
An absolute-error criterion chooses the subproblem error a priori \citep{eckstein1992douglas}, 
while a relative error criterion requires the subproblem error 
to decrease as the algorithm nears convergence, for example, 
by setting the error target at each iteration proportional to 
$\|u^{k+1} - u^k - \rho (z^{k+1} - z^k)\|$ \citep{eckstein2018relative}. 

\paragraph{Approximations used by GeNI-ADMM.}
The GeNI-ADMM framework allows for any combination of the three strategies: 
augmented Lagrangian linearization, function approximation, and inexact subproblem solve.
Consider the generalized second-order approximation to $f$
\begin{align}
\label{eq:gen-newt-appx}
f(x) \approx f_1(x) + f_2(\tilde x^k)+\langle \nabla f_2(\tilde x^k ), x - \tilde x^k \rangle + \frac{1}{2}\|x - \tilde x^k\|^2_ {\Theta^k},
\end{align}
where $f_1+f_2 = f$ and $\{\Theta^k\}_{k\ge 1}$ is a sequence of psd matrices 
that approximate the Hessian of $f$.
GeNI-ADMM uses this approximation in the $x$-subproblem,
resulting in the new subproblem
\begin{equation}
\label{eq:gen-newt-step}
\begin{aligned}
\tilde x^{k+1} = \underset{x\in X}{\textrm{argmin}} & \{ 
f_1(x) + \langle\nabla f_2(\tilde x^k ),x - \tilde x^k\rangle + \frac{1}{2}\|x - \tilde x^k\|^2_ {\Theta^k} \\
& +\frac{\rho}{2}\|Mx +N\tilde z^k + \tilde u^k\|^2 \}.
\end{aligned}
\end{equation}
 
We refer to \eqref{eq:gen-newt-step} as a generalized Newton step. 
The intuition for the name is made plain when $X = \mathbb{R}^d$, $f_1 = 0, f_2 = f$
in which case the update becomes
\begin{align}
\label{eq:gen-newt-update}
&\left(\Theta^k+\rho M^{T}M\right)\tilde{x}^{k+1} = \Theta^k \tilde x^k-\nabla f(\tilde x^k)
-\rho M^{T}(\tilde u^k - \tilde z^k) .
\end{align}
Equation \eqref{eq:gen-newt-update} shows that the $x$-subproblem reduces to a linear system solve, just like the Newton update.
We can also interpret GeNI-ADMM as a linearized proximal augmented Lagrangian (P-ALM) method \citep{hermans2019qpalm,hermans2022qpalm}. 
From this point of view, GeNI-ADMM replaces $f$ in the P-ALM step by its linearization and adds a specialized penalty defined by the $\Theta^k$-norm.

GeNI-ADMM also replaces the $z$-subproblem of ADMM with the following problem
\begin{equation}
\label{eq:prox-z-subprob}
    \tilde{z}^{k+1} = \underset{z\in Z}{\argmin} \{g(z)+ \frac{\rho}{2}\|M\tilde x^{k+1} +Nz + \tilde u^k\|^2+\frac{1}{2}\|z-\tilde z^k\|^2_{\Psi^k}\}.
\end{equation}
Incorporating the quadratic $\Psi^k$ term allows us to linearize the augmented Lagrangian term, which is useful when $N$ is very complicated.

However, even with the allowed approximations, it is unreasonable to assume that \eqref{eq:gen-newt-step},~\eqref{eq:prox-z-subprob} are solved exactly at each iteration. 
Indeed, the hallmark of the NysADMM scheme from \cite{pmlr-v162-zhao22a} is that it solves \eqref{eq:gen-newt-step} inexactly (but efficiently) using a randomized linear system solver.
Thus, the GeNI-ADMM framework allows for inexactness in the $x$ and $z$-subproblems, as seen in \cref{alg-gen_newton_admm}. 
\begin{algorithm}[H]
	\centering
	\caption{Generalized Newton Inexact ADMM (GeNI-ADMM)\label{alg-gen_newton_admm}}
	\begin{algorithmic}
		\Require{penalty parameter $\rho$, sequence of psd matrices $\{\Theta^k\}_{k \ge 1}, \{\Psi^k\}_{k\ge 1}$, forcing sequences $\{\varepsilon_x^k\}_{k\geq 1}, \{\varepsilon_z^k\}_{k\geq 1},$}
		\Repeat
		\State{$\tilde{x}^{k+1} \overset{\varepsilon_x^k}{\approx}\underset{x \in X}{\argmin}\{ f_1(x) + \langle \nabla f_2(\tilde x^k ), x - \tilde x^k\rangle+\frac{1}{2}\|x - \tilde x^k\|^2_ {\Theta^k} + \frac{\rho}{2}\|Mx+N\tilde z^k+ \tilde u^k\|^2\}$}
		\State{$\tilde{z}^{k+1} \overset{\varepsilon_z^k}{\approxeq} \underset{z\in Z}{\argmin} \{g(z)+ \frac{\rho}{2}\|M\tilde x^{k+1} +Nz + \tilde u^k\|^2+\frac{1}{2}\|z-\tilde z^k\|^2_{\Psi^k}\}$}
		\State{$\tilde u^{k+1} = \tilde u^k + M\tilde x^{k+1} +N\tilde z^{k+1}$}
		\Until{convergence}
		\Ensure{solution $(x^\star, z^\star)$ of problem \eqref{eq:ADMMProb}}
	\end{algorithmic}
\end{algorithm}
\Cref{alg-gen_newton_admm} differs from ADMM (\cref{alg-admm}) in that a) the $x$ and $z$-subproblems are now given by \eqref{eq:gen-newt-step} and \eqref{eq:prox-z-subprob} and b) both subproblems may be solved inexactly.
Given the inexactness and the use of the generalized Newton step in place of the original $x$-subproblem, we refer to \cref{alg-gen_newton_admm} as Generalized Newton Inexact ADMM (GeNI-ADMM). 
The inexactness schedule is controlled by the \textit{forcing sequences} $\{\varepsilon_x^k\}_{k\geq 1}, \{\varepsilon_z^k\}_{k\geq 1}$  that specify how accurately the $x$ and $z$-subproblems are solved at each iteration.
These subproblems have different structures that require different notions of accuracy. To distinguish between them, we make the following definition. 

\begin{definition}[$\varepsilon$-minimizer and $\varepsilon$-minimum]
\label{def:TypeInexactness}
Let $h:\mathcal {T}\mapsto \mathbb{R}$ be strongly-convex and let $t^\star = \argmin_{t'\in \mathcal{T}} h(t')$. 
\begin{itemize}
    \item ($\varepsilon$-minimizer) Given $t \in \mathcal{T}$, we say $t$ is an $\varepsilon$-\emph{minimizer} of $\textrm{minimize}_{t\in \mathcal{T}}h(t)$ and write
    \[t \epsappx \argmin_{t'\in \mathcal{T}} h(t') \quad \textrm{if and only if} \quad \|t-t^\star\| \leq \varepsilon.\]
    In words, $t$ is nearly equal to the argmin of $h(t)$ in set $\mathcal T$.

     \item ($\varepsilon$-minimum) Given $t \in \mathcal{T}$, we say $t$ gives an $\varepsilon$-\emph{minimum} of $\textrm{minimize}_{t\in \mathcal{T}}h(t)$ and write
    \[t \overset{\varepsilon}{\approxeq} \argmin_{t'\in \mathcal{T}} h(t') \quad \textrm{if and only if} \quad h(t)-h(t^\star) \leq \varepsilon.\]
    In words, $t$ produces nearly the same objective value as the argmin of $h(t)$ in set $\mathcal T$.
    \end{itemize}
\end{definition}
Thus, from \cref{def:TypeInexactness} and \cref{alg-gen_newton_admm}, we see for each iteration $k$ that $\tilde{x}^{k+1}$ is an $\varepsilon_x^k$-minimizer of the $x$-subproblem, while $\tilde{z}^{k+1}$ gives an $\varepsilon_z^k$-minimum of the $z$-subproblem. 

\subsection{Related work}
The literature on the convergence of ADMM and its variants is vast, 
so we focus on prior work most relevant to our setting.
\cref{table:relatedwork} lists some of the prior work that developed and analyzed the approximation strategies described in \cref{section:GeneralizedNewtADMM}. 
GeNI-ADMM differs from all prior work in \cref{table:relatedwork} by allowing (almost) all these approximations and more.
It also provides explicit rates of convergence to support choices between algorithms.
Moreover, many of these algorithms can be recovered from GeNI-ADMM.

\begin{table}[htpb] 
\scriptsize
\centering 
    \begin{tabular}{C{2cm}C{2cm}C{2cm}C{2cm}C{2cm}C{2cm}}
        \textbf{Reference} & \textbf{Convergence \newline rate} & \textbf{Problem \newline class} & \textbf{Augmented Lagrangian Linearization}  & \textbf{Function \newline approximation}  & \textbf{Subproblem \newline inexactness}  \\ \hline
        \citet{eckstein1992douglas}  & \redX & Convex & \redX & \redX & $x,z$ \newline absolute error \\ \hline
        \citet{he20121} & Ergodic $\bigO(1/t)$  & Convex & $x$ &\redX & \redX\\ \hline
        \citet{monteiro2013iteration}  & Ergodic $\bigO(1/t)$ & Convex & \redX & \redX & \redX \\ \hline
        \citet{ouyang2013stochastic}  & $\bigO(1/\sqrt{t})$ & Stochastic \newline convex & \redX & $f$ \newline stochastic first-order & \redX \\ \hline
        \citet{ouyang2015accelerated}  & Ergodic $\bigO(1/t)$ & Convex & $x$ & $f$ \newline first-order & \redX \\ \hline
        \citet{deng2016global} & Linear & Strongly convex & $x,z$ & \redX & \redX \\ \hline
        \citet{eckstein2018relative} & \redX & Convex & \redX & \redX & $x$ \newline relative error \\ \hline
        \citet{hager2020convergence} &  Ergodic $\bigO(1/t),$ \newline $\bigO(1/t^2)$ & Convex, \newline Strongly convex & $x$  & \redX & $x$ \newline relative error \\ \hline
        \citet{yuan2020discerning} &  Locally \newline linear & Convex & $x,z$  & \redX & \redX \\ \hline
        \citet{pmlr-v162-zhao22a} & \redX \newline Convergence \newline for quadratic $f$ only &Convex & \redX & $x$ & $x$ absolute error \newline \\ \hline
        \citet{ryu2022large} & \redX & Convex & $x,z$ & $f,g$ \newline partial first-order & \redX \\ \hline
        \textbf{This work} & Ergodic $\bigO(1/t),$ \newline Linear & Convex, \newline Strongly convex & $x,z $& $f$ \newline partial generalized second-order & $x, z$ \newline absolute error \\ \hline
    \end{tabular}
    \caption{A structured comparison of related work on the convergence of ADMM and its variants. The ``$x$'' (``$z$'') in the table denotes that a paper uses the corresponding strategy of the column to simplify the $x$-subproblem ($z$-subproblem). In the ``\textbf{Function approximation}'' column, the ``$f$'' (``$g$'') indicates that a paper approximates function ``$f$'' (``$g$'') in the $x$-subproblem ($z$-subproblem). ``stochastic first-order'' means a paper uses first-order function approximation but replace the gradient term with a stochastic gradient. ``partial generalized second-order'' means a paper uses second-order function approximation as \eqref{eq:gen-newt-step}, but $H$ is not necessarily the Hessian.} 
    \label{table:relatedwork}
\end{table}

\subsubsection{Algorithms recovered from GeNI-ADMM}
\label{subsection:recoveredadmm}
Various ADMM schemes in the literature can be recovered by appropriately selecting the parameters in \cref{alg-gen_newton_admm}. 
Let us consider a few special cases to highlight important prior work on approximate ADMM, 
and provide concrete intuition for the general framework provided by \cref{alg-gen_newton_admm}.

\paragraph{NysADMM}
The NysADMM scheme \citep{pmlr-v162-zhao22a} assumes the problem is unconstrained $M=I, N = -I$.
\citet{diamandis2023genios} have extended NysADMM to the constrained setting with the GeNIOS solver, which we will discuss more in \cref{subsection:NysADMM}.
For simplicity of exposition, we focus on the original NysADMM for unconstrained problems.  
GeNI-ADMM specializes to NysADMM taking $\Theta^k = \eta (H_f^k+\sigma I)$, and $\Psi^k = 0$, 
where $H_f^k$ denotes the Hessian of $f$ at the $k$th iteration.
Unlike the original NysADMM scheme of \cite{pmlr-v162-zhao22a}, we include the regularization term $\sigma I$, where $\sigma>0$.
This inclusion is required for theoretical analysis but seems unnecessary in practice
(see \cref{subsection:NysADMM} for a detailed discussion). 
Substituting this information into \eqref{eq:gen-newt-update} leads to the following update for the $x$-subproblem. 
\begin{equation}
\label{eq:nysadmm-update}
\tilde{x}^{k+1} = \tilde x^k-\left(\eta H_f^k+(\rho +\eta \sigma)I\right)^{-1}\left(\nabla f(\tilde x^k)+\rho ( \tilde x^k-\tilde z^k+\tilde u^k)\right).
\end{equation}

\paragraph{Sketch-and-solve ADMM}
If $M = I, N=-I$, $\Psi^k = 0$, 
and $\Theta^k$ is chosen to be a matrix such that $\Theta^k+\rho I$ is easy to factor, 
we call the resulting scheme \emph{sketch-and-solve ADMM}. 
The name \emph{sketch-and-solve ADMM} is motivated by the fact that such a $\Theta^k$ can often be obtained via sketching techniques.
However, the method works for any $\Theta^k$, not only ones constructed by sketching methods.
The update is given by
\begin{equation*}
\tilde{x}^{k+1} = \tilde x^k-\left(\Theta^k+\rho I\right)^{-1}\left(\nabla f(\tilde x^k)+\rho (\tilde x^k-\tilde z^k+\tilde u^k)\right).
\end{equation*}
We provide the full details of sketch-and-solve ADMM in \cref{subsection:SketchAndSolveADMM}.
To our knowledge, sketch-and-solve ADMM has not been previously proposed or analyzed.

\paragraph{Proximal/Generalized ADMM}
Set $f_1 = f, f_2 =0$, and fix $\Theta^k = \Theta$ and $\Psi^k = \Psi$, where $\Theta$ and $\Psi$ are symmetric positive definite matrices. Then GeNI-ADMM simplifies to
\begin{align*}
    &\tilde x^{k+1} = \argmin_{x\in X}\left\{f(x)+\frac{\rho}{2}\|Mx+N\tilde z^k+\tilde u^k\|^2+\frac{1}{2}\|x-\tilde x^k\|_{\Theta}^2\right\},\\
    &\tilde z^{k+1} = \argmin_{z\in Z}\left\{g(z)+\frac{\rho}{2}\|M\tilde x^{k+1}+Nz+\tilde u^k\|^2+\frac{1}{2}\|z-\tilde z^k\|_{\Psi}^2\right\},
\end{align*}
which is the Proximal/Generalized ADMM of \citet{deng2016global}.

\paragraph{Linearized ADMM}
Set $\Theta^k = \alpha^{-1}I-\rho M^TM$ and $\Psi = \beta^{-1}I-\rho N^{T}N$, for all $k\geq 1$. Then GeNI-ADMM simplifies to
\begin{align*}
&\tilde x^{k+1} = \prox_{\alpha f}\left(\tilde x^k-\alpha\rho M^{T}(M\tilde x^k+N\tilde z^k+\tilde u^k)\right), \\
&\tilde z^{k+1} = \prox_{\beta g}\left(\tilde z^k-\beta\rho N^{T}(M\tilde x^{k+1}+N\tilde z^k+\tilde u^k)\right),
\end{align*}
which is exactly Linearized ADMM \citep[\S 4.4.2]{parikh2014proximal}.

\paragraph{Primal Dual Hybrid Gradient}
The Primal Dual Hybrid Gradient (PDHG) or Chambolle-Pock algorithm  of \citet{chambolle2011first} is a special case of GeNI-ADMM since PDHG is a special case of Linearized ADMM (see section 3.5 of \cite{ryu2022large}).

\paragraph{Gradient descent ADMM}
Consider two linearization schemes from \citet{ouyang2015accelerated}.
The first scheme we  call \emph{gradient descent ADMM} (GD-ADMM) is useful when it is cheap to solve a least-squares problem with $M$.
GD-ADMM is obtained from GeNI-ADMM by setting $f_1 = 0, f_2 = f$, and $\Theta^k = \eta I$ for all $k$. 
The $\tilde{x}^{k+1}$ update \eqref{eq:gen-newt-update} for GD-ADMM simplifies to 
\begin{equation}
    \label{eq:gd-admm-update}
    \tilde x^{k+1} = \tilde x^k-(\rho M^{T}M+\eta I)^{-1}\left(\nabla f(\tilde x^k)+\rho M^T(M\tilde x^k+N\tilde z^k+\tilde u^k)\right).
\end{equation}
The second scheme, \emph{linearized-gradient descent ADMM} (LGD-ADMM), is useful when $M$ is not simple, so that the update \eqref{eq:gd-admm-update} is no longer cheap.
To make the $x$-subproblem update cheaper, it linearizes the augmented Lagrangian term by setting $\Theta^k = \alpha^{-1}I-\rho M^{T}M$ for all $k$ in addition to linearizing $f$. In this case, \eqref{eq:gen-newt-update} yields the $\tilde{x}^{k+1}$ update 
\begin{equation}
\label{eq:pgd-admm-update}
    \tilde x^{k+1} = \tilde x^k-\alpha\left(\nabla f(\tilde x^k)+\rho M^{T}(M\tilde x^k+N\tilde z^k+\tilde u^k)\right).
\end{equation}
Observe in the unconstrained case, when $M = I$,  the updates \eqref{eq:gd-admm-update} and \eqref{eq:pgd-admm-update} are equivalent and generate the same iterate sequences when initialized at the same point \citep{zhao2021automatic}. Indeed, they are both generated by performing a gradient step   
on the augmented Lagrangian \eqref{eq:gen-newt-step}, for suitable choices of the parameters.
Notably, this terminology differs from \cite{ouyang2015accelerated}, who refer to \eqref{eq:gd-admm-update} as ``linearized ADMM'' (L-ADMM) and \eqref{eq:pgd-admm-update} as ``linearized preconditioned ADMM'' (LP-ADMM). 
We choose our terminology to emphasize that GD-ADMM accesses $f$ via its gradient, as in the literature the term ``linearized ADMM'' is usually reserved for methods that access $f$ through its prox operator \citep{he20121,parikh2014proximal,deng2016global}.

In the remainder of this paper, we will prove convergence of \cref{alg-gen_newton_admm} under appropriate hypotheses on the sequences $\{\Theta^k\}_{k\geq 1}, \{\Psi^k\}_{k\geq 1}, \{\varepsilon_x^k\}_{k\geq 1},$~and~$\{\varepsilon_z^k\}_{k\geq 1}$.

\section{Technical preliminaries and assumptions}
\label{section:TechPrelimsAssump}
We introduce some important concepts that will be central to our analysis.
The first is $\Theta$-\textit{relative smoothness}, which is crucial to establish that GeNI-ADMM benefits from curvature information provided by the Hessian.
\begin{definition}[$\Theta$-relative smoothness]
\label{def:rel-smoothness}
We say $f: \mathcal D \to \reals$ is $\Theta$-{relatively smooth} with respect to 
the bounded function $\Theta: \mathcal D \to \sym_+^n$ 
if there exists $\hat{L}_\Theta > 0$ such that for all $x,y \in \mathcal{D}$
\begin{equation}
\label{eq:rel-smoothness}
    f(x)\leq f(y)+\langle \nabla f(y),x-y\rangle+\frac{\hat{L}_\Theta}{2}\|x-y\|^2_{\Theta(y)}.
\end{equation}
\end{definition}
That is, the function $f$ is smooth with respect to the $\Theta$-norm.
\Cref{def:rel-smoothness} generalizes \emph{relative smoothness}, 
introduced in \cite{gower2019rsn} to analyze Newton's method.
The definition in \cite{gower2019rsn} takes $\Theta$ to be the Hessian of $f$, $H_f$. 
When $f$ belongs to the popular family of generalized linear models,
then~\eqref{eq:rel-smoothness} holds with a value of $\hat{L}_{H_f}$ 
that is independent of the conditioning of the problem~\citep{gower2019rsn}.
For instance, if $f$ is quadratic and $\Theta(y) = H_f$, 
then \eqref{eq:rel-smoothness} holds with equality for $\hat{L}_{H_f} = 1$. 
Conversely, if we take $\Theta = I$, which corresponds to GD-ADMM, 
then $\hat{L}_\Theta = L_f$, the smoothness constant of $f$. 
Our theory uses the fact that $\hat{L}_{\Theta}$ is much smaller than the smoothness constant $L_f$
for methods that take advantage of curvature,
and relies on $\hat L_\Theta$ to characterize the faster convergence speed of these methods.

The other important idea  we need is the notion of an \textit{$\varepsilon$-subgradient} \citep{bertsekas2003convex,hiriart1993convex}.
\begin{definition}[$\varepsilon$-subgradient]
\label{def:epssubgrad}
 Let $r: \mathcal D \to \reals $ be a convex function and $\varepsilon > 0$. 
 We say that $s \in \mathcal{D}^*$ is an $\varepsilon$-subgradient for $r$ at $z \in \mathcal D$
 if, for every $z' \in \mathcal D$, we have
\[r(z')-r(z)\geq \langle s, z'-z\rangle -\varepsilon.\]
We denote the set of $\varepsilon$-subgradients for $r$ at $z$ by $\partial_{\varepsilon} r(z)$.
\end{definition}
Clearly, any subgradient is an $\varepsilon$-subgradient, so \cref{def:epssubgrad} provides a natural weakening of a subgradient.
The $\varepsilon$-subgradient is critical for analyzing $z$-subproblem inexactness, and our usage in this context is inspired by the convergence analysis of inexact proximal gradient methods \citep{schmidt2011convergence}.
We shall need the following proposition whose proof may be found in~\cite[Proposition 4.3.1]{bertsekas2003convex}.
\begin{proposition}
\label{proposition:epssubgrad-properties}
Let $r$, $r_1$, and $r_2$ be convex functions. Then for any $z$, the following holds:
\begin{enumerate}
    \item $0\in \partial_\varepsilon r(z)$ if and only if $r(z) \overset{\varepsilon}{\approxeq}\argmin_{z'} r(z')$, that is $z$ gives an $\varepsilon$-minimum of $\textup{minimize}_{z'} r(z')$.
    \item $ \partial_\varepsilon(r_1+r_2)(z) \subset \partial_\varepsilon r_1(z)+ \partial_\varepsilon r_2(z).$
\end{enumerate}
\end{proposition}
With \cref{proposition:epssubgrad-properties} recorded, 
\ifpreprint
we prove the following lemma in \cref{subsection:SubLinConvPfs} of the appendix,
\fi
which will play a critical role in establishing the convergence of GeNI-ADMM.  

\begin{lemma}
\label{lemma:zsubprob}
Let $\tilde{z}^{k+1}$ give an $\varepsilon^k_z$-minimum of the $z$-subproblem. 
Then there exists an $\tilde s$ with $\|\tilde s\|\leq C\sqrt{\varepsilon_z^k}$
such that
\[
-\rho N^{T}\tilde u^{k+1}+\Psi_k(\tilde z^k-\tilde z^{k+1})+\tilde s \in \partial_{\varepsilon_z^k}g(\tilde{z}^{k+1}).
\]
\end{lemma}

\subsection{Assumptions}
In this section, we present the main assumptions required by our analysis.
\begin{assumption}[Existence of saddle point]
\label{assp:existence}
There exists an optimal primal solution $(x^\star, z^\star) \in X \times Z $ for \eqref{eq:ADMMProb} 
and an optimal dual solution $u^\star \in U$ such that $(x^\star, z^\star, u^\star)$ is a saddle point of \eqref{eq:saddlepoint}. 
Here, $U \subset \mathcal Z$ is a closed convex set and $\rho U = Y$. 
We denote the optimal objective value of \eqref{eq:ADMMProb} as $p^\star$.
\end{assumption}
Assumption \ref{assp:existence} is standard and merely assumes that \eqref{eq:ADMMProb} has a solution.

\begin{assumption}[Regularity of $f$ and $g$] 
\label{assp:regularity}
The function $f$ is  twice-continuously differentiable and is $1$-relatively smooth with respect to $\Theta$. 
The function $g$ is finite-valued, convex, and lower semi-continuous. 
\end{assumption}
Assumption \ref{assp:regularity} is also standard. It ensures that it makes sense to talk about the Hessian of $f$ and that $f$ is relatively smooth.
Note assuming $\hat L_\Theta = 1$ is without loss of generality, 
for we can always redefine $\Theta' = \hat L_\Theta \Theta$, 
and $f$ will be $1$-relatively smooth with respect to $\Theta'$.  
\begin{assumption}[Forcing sequence summability and approximate subproblem oracles]
\label[assumption]{assp:x-subprob-orc}
Let $\{\varepsilon^k_x\}_{k\ge 1}$ and $\{\varepsilon^k_z\}_{k \ge 1}$ be given forcing sequences, we assume they satisfy 
\[\mathcal{E}_x = \sum_{k=1}^\infty\varepsilon_x^k<\infty, \quad \mathcal{E}_z = \sum^\infty_{k=1} \sqrt{\varepsilon_z^k}<\infty.\]
Further we define the constants $K_{\varepsilon_x} \coloneqq \sup_{k\geq 1}\varepsilon_x^k$, $K_{\varepsilon_z} \coloneqq \sup_{k\geq 1}\sqrt{\varepsilon_z^k}$. 
Observe $K_{\varepsilon_x}$ and $K_{\varepsilon_z}$ are finite owing to the summability hypotheses. 
Moreover, we assume \cref{alg-gen_newton_admm} is equipped with oracles for solving the $x$ and $z$-subproblems, which at each iteration produce approximate solutions $\tilde{x}^{k+1}$, $\tilde z^{k+1}$ satisfying:
\begin{align*}
    & \tilde x^{k+1} \xepsappx  \underset{x \in X}{\argmin}\{ f_1(x) + \langle \nabla f_2(\tilde x^k ), x - \tilde x^k\rangle+\frac{1}{2}\|x - \tilde x^k\|^2_ {\Theta^k} + \frac{\rho}{2}\|Mx+N\tilde z^k + \tilde u^k\|^2\},\\
    & \tilde z^{k+1} \zfepsappx \underset{z\in Z}{\argmin} \{g(z)+ \frac{\rho}{2}\|M\tilde x^{k+1} +Nz + \tilde u^k\|^2+\frac{1}{2}\|z-\tilde z^k\|^2_{\Psi^k}\}. \end{align*}
\end{assumption}
The conditions on the $x$ and $z$ subproblem oracles are consistent with those of \cite{eckstein1992douglas}, 
which requires the sum of the errors in the subproblems to be summable.
The approximate solution criteria of Assumption \ref{assp:x-subprob-orc} are also easily met in practical applications, as the subproblems are either simple enough to solve exactly, or can be efficiently solved approximately via iterative algorithms.
For instance, with LGD-ADMM the $x$-subproblem is simple to solve enough exactly, 
while for NysADMM the $x$-subproblem is efficiently solved via conjugate gradient with a randomized preconditioner. 

\begin{assumption}[Regularity of $\{\Theta^k\}_{k \ge 1}$ and $\{\Psi^k\}_{k\ge 1}$]
\label{assp:metric_evol}
We assume there exists constants $\theMax,~\theMin,~\psiMax,$ and $\nu>0$, such that
\begin{equation}
\label{eq:mat_reg_assmp}
\theMax I\succeq \Theta^k\succeq \theMin I,~\psiMax I\succeq \Psi^k, ~\Psi^k+\rho N^TN\succeq \nu I, \quad \text{for all}~k\geq 1.
\end{equation}
Moreover, we also assume that The sequences $\{\Theta^k\}_{k \ge 1}, \{\Psi^k\}_{k\ge 1}$ satisfy 
\begin{align}
	\label{eq:Hcond}
    &\|\tilde x^k-x^\star\|^2_{\Theta^k}\leq (1+\zeta^{k-1})\|\tilde x^k-x^\star\|^2_{\Theta^{k-1}}, \\
    &\|\tilde z^k-z^\star\|^2_{\Psi^k}\leq (1+\zeta^{k-1})\|\tilde z^k-z^\star\|^2_{\Psi^{k-1}}, \nonumber
\end{align}
where $\{\zeta^k\}_{k \ge 1}$ is a non-negative summable sequence, that is,  $\sum^{\infty}_{k=1}\zeta^k <\infty$. 
\end{assumption}
The first half of Assumption \ref{assp:metric_evol} is standard, and essentially requires that $\Theta^k$ and $\Psi^k+\nu N^TN$ define norms. 
For most common choices of $\Theta^k$ and $\Psi^k$ these assumptions are satisfied, they are also readily enforced by adding on a small multiple of the identity.
The addition of a small regularization term is common practice in popular ADMM solvers, such as OSQP \citep{osqp}, as it avoids issues with degeneracy and increases numerical stability. 

The second part of the assumption requires that $\Theta^{k}$ ($\Psi^k)$ and $\Theta^{k-1}$ ($\Psi^{k-1}$) eventually not differ much on the distance of $\tilde x^k$ ($\tilde z^k$) to $x^\star$ ($z^\star$).
Assumptions of this form are common in analyses of optimization algorithms that use variable metrics. 
For instance, \citet{he2002new} which develops an alternating directions method for monotone variational inequalities,
assumes their equivalents of the sequences $\{\Theta^k\}$ and $\{\Psi^k\}$ satisfy
\begin{equation}
\label{eq:metric}
    (1-\zeta^{k-1})\Theta^{k-1}\preceq \Theta^k \preceq (1+\zeta^{k-1})\Theta^{k-1},~(1-\zeta^{k-1})\Psi^{k-1}\preceq \Psi^k \preceq (1+\zeta^{k-1})\Psi^{k-1}.
\end{equation}
More recently, \citet{rockafellar2023generic} analyzed the proximal point method with variable metrics, under the assumption the variable metrics satisfy \eqref{eq:metric}. 
Thus, Assumption \ref{assp:metric_evol} is consistent with the literature and is, in fact, weaker than prior work, as \eqref{eq:metric} implies \eqref{eq:Hcond}.
Assumption \ref{assp:metric_evol} may always be enforced by changing the $\Theta^k$'s and $\Psi^k$'s only finitely many times.

We also define the following constants which shall be useful in our analysis,
\[
\tau_\zeta \coloneqq \prod_{k\geq 2}(1+\zeta^k), \quad \mathcal{E}_{\zeta} = \tau_\zeta\left(\sum^{\infty}_{k=1}\zeta^k \right)<\infty.
\]
Note Assumption \ref{assp:metric_evol} implies $\tau_\zeta <\infty $.
Moreover, for any $k\geq 1$  it holds that 
\[
\|\tilde x^k-x^\star\|^2_{\Theta^k}\leq\tau_\zeta\|\tilde x^k-x^\star\|^2_{\Theta^1},~\|\tilde z^k-z^\star\|^2_{\Psi^k}\leq \tau_\zeta \|\tilde z^k-z^\star\|^2_{\Psi^1}.
\] 
\section{Sublinear convergence of GeNI-ADMM}
\label{section:SubLinConvergence}
This section establishes our main theorem, \cref{thm:generalconv}, 
which shows that \cref{alg-gen_newton_admm} enjoys the same $\mathcal O(1/t)$-convergence rate as standard ADMM. 
\begin{theorem}[Ergodic convergence]
    \label{thm:generalconv} 
    Define constants  
    $d_{x^\star,\Theta^1} = \|\tilde x^1-x^\star\|_{\Theta^1}$, $d_{u^\star} = \|\tilde u^1 - u^\star\|$, $d_{z^\star,\Psi^1_{\rho, N}} = \|\tilde z^1 - z^\star\|_{\Psi^1+\rho N^TN}$.     
    Let $p^\star$ denote the optimum of \eqref{eq:ADMMProb}. 
    For each $t\geq 1$, denote $\bar x^{t+1} = \frac{1}{t}\sum^{t+1}_{k=2}\tilde x^k$, and 
    $\bar z^{t+1} = \frac{1}{t}\sum^{t+1}_{k=2}\tilde z^k$, where $\{\tilde x^k\}_{k\geq 1}$ and $\{\tilde z^k\}_{k\geq 1}$ are the iterates produced by \cref{alg-gen_newton_admm} with forcing sequences $\{\varepsilon^k_x\}_{k \ge 1}$ and $\{\varepsilon^k_z\}_{k \ge 1}$ with $\tilde u^1 = 0$ and $M\tilde x^1=-N\tilde z^1 $.  Instate Assumptions \ref{assp:existence}-\ref{assp:metric_evol}. Then, the suboptimality gap satisfies
    \small
    \begin{equation*}
	f(\bar x^{t+1})+g(\bar z^{t+1})-p^{\star}\leq 
		\frac{1}{t}\left(\frac{1}{2}d^2_{x^\star,\Theta^1}+\frac{1}{2}d_{z^\star,\Psi^1_{\rho, N}}^{2}+ C_x\mathcal{E}_x+C_z\mathcal{E}_z+C_\zeta\mathcal{E}_{\zeta}\right) =: \frac 1 t \Gamma,\\
	\end{equation*}
    \normalsize
    Furthermore, the feasibility gap satisfies 
    \small
    \begin{equation*}
        \|M\bar x^{t+1}+N\bar z^{t+1}\| \le \frac{2}{t}\sqrt{\frac{\Gamma}{\rho}+d^2_{u^\star}},
    \end{equation*}
    Consequently, after $\mathcal O(1/\epsilon)$ iterations,
	\begin{equation*}
	f(\bar x^{t+1})+g(\bar z^{t+1})-p^{\star}\leq \epsilon, \; \text{and} \;\;\; \|M \bar x^{t+1} +N\bar z^{t+1}\| \le \epsilon. 
	\end{equation*}
\end{theorem}
\Cref{thm:generalconv} shows that, with a constant value of $\eta$ and appropriate forcing sequences, the suboptimality gap and the feasibility residuals both go to zero at a rate of $\mathcal O(1/t)$. 
Hence, the overall convergence rate of ADMM is preserved despite all the approximations involved in GeNI-ADMM.
 
\subsection{The benefits of using a better approximation}
To see this difference, consider the case when $f$ is a convex quadratic function and $M = I$. 
Let $H_f$ be the Hessian of $f$. 
Consider (a) NysADMM $\Theta^k = H_f$ \eqref{eq:nysadmm-update} and (b) GD-ADMM $\Theta^k = I$ \eqref{eq:gd-admm-update}.
Further, suppose that NysADMM and GD-ADMM are initialized at $0$.
The rates of convergence guaranteed by \cref{thm:generalconv} for \cref{alg-gen_newton_admm} for both methods are outlined in \cref{table:convergerates}. 
\begin{table}[tbhp]
\caption{Convergence rate comparison of NysADMM and GD-ADMM when initialized at $0$ for quadratic $f$.}
\label{table:convergerates}
\vskip 0.15in
\begin{small}
\begin{tabular}{lcc}
    \hline
     Method & NysADMM $\left(\Theta^k = H_f\right)$ & GD-ADMM $\left(\Theta^k = I\right)$  \\
    \hline
    Feasibility gap & $\mathcal O\left(\frac{1}{t}\left(\sqrt{\frac{2}{\rho}\|x^\star\|^2_{H_f}}\right)\right)$ & $\mathcal O\left(\frac{1}{t}\left(\sqrt{\frac{2}{\rho}\lambda_1\left(H_f\right)\|x^\star\|^2}\right)\right)$ \\
    Suboptimality gap & $\mathcal O\left(\frac{1}{2t}\|x^\star\|^2_{H_f}\right)$ & $\mathcal O\left(\frac{1}{2t}\lambda_1\left(H_f\right)\|x^\star\|^2\right)$ \\
    \hline
\end{tabular}
\end{small}
\vskip -0.1in
\end{table}
In the first case, the relative smoothness constant satisfies $\hat L_\Theta = 1$, while in the second, $\hat L_\Theta = \lambda_1\left(H_f\right)$, which is the largest eigenvalue of $H_f$. 
Comparing the rates in \cref{table:convergerates}, we see NysADMM improves over GD-ADMM whenever 
\[
\lambda_1\left(H_f\right)\|x^\star\|^2 \ge \|x^\star\|^2_{H_f}.
\]
Hence NysADMM improves significantly over GD-ADMM when $H_f$ exhibits a decaying spectrum, 
provided $x^\star$ is not concentrated on the top eigenvector of $H_f$.
We formalize the latter property in the following definition.

\begin{definition}[$\mu_\mathbf{V}$-incoherence]
\label{def:muVIncoho}
Let $\mathbf V$ be the eigenbasis of $H_f$. 
We say $x^\star$ is $\mu_\mathbf{V}$-incoherent if there exists $1 \leq \mu_\mathbf{V} < d$ such that:
\begin{equation}
\label{eq:IncohAssmp}
\sup_{1\leq i \leq d}|\langle \mathbf{v}_i, x^\star \rangle |^2 \leq \mu_\mathbf{V} \frac{\|x^\star\|^2}{d}.
\end{equation}
\end{definition}
\Cref{def:muVIncoho} is a weak form of the incoherence condition from compressed sensing and matrix completion, which plays a key role in signal  and low-rank matrix recovery \citep{candes2007sparsity,candes2012exact}. 
In words, $x^\star$ is $\mu_\mathbf{v}$-incoherent if its energy is not solely concentrated on the top eigenvector of $H_f$ and can be expected to hold generically.
 The parameter $\mu_\mathbf{v}$ controls the allowable concentration.  
When $\mu_\mathbf{v} = 0$, $x^\star$ is orthogonal to $\mathbf{v}$, so its energy is distributed amongst the other eigenvectors.
Conversely, the closer $\mu_\mathbf{v}$ is to $1$, the more $x^\star$ is allowed to concentrate on $\mathbf{v}$.

Using $\mu_\mathbf{v}$-incoherence, we can say more about how NysADMM improves on GD-ADMM.
\begin{proposition}
\label{prop:mu-incoherence}
Suppose $x^{\star}$ is $\mu_\mathbf{v}$-incoherent and $\sigma = \tau \lambda_1(H_f)$ where $\tau \in (0,1)$. Then, the following bound holds
\[\frac{\lambda_{1}(H_f)\|x^\star\|^2}{\|x^\star\|^2_{H_f}} \geq 
\frac{d}{\mu_{\mathbf{V}}\textup{intdim}(H_f)}\]
Hence if $\mu_\mathbf{v}+\lambda_2(H_f)/\lambda_1(H_f)\leq (1-\tau)\alpha^{-1}$, 
where $\alpha \geq 1$, then
\[\frac{\lambda_{1}(H_f)\|x^\star\|^2}{\|x^\star\|^2_{H_f}} \geq \alpha.\]
\end{proposition}

\ifpreprint
The proof of \cref{prop:mu-incoherence} is given in \cref{subsection:SubLinConvPfs}.
\fi
\cref{prop:mu-incoherence} shows when $x^\star$ is $\mu_\mathbf{v}$-incoherent and $H_f$ has a decaying spectrum, NysADMM yields a significant improvement over GD-ADMM. 
As a concrete example, consider when $\alpha = 2$, 
then \Cref{prop:mu-incoherence} implies the ergodic convergence of NysADMM is twice as fast as GD-ADMM. 
We observe this performance improvement in practice; see \cref{section:experiments} for corroborating numerical evidence.
Just as Newton's method improves on gradient descent for ill-conditioned problems, 
NysADMM is less sensitive to ill-conditioning than GD-ADMM.

\subsection{Our approach}
To prove \cref{thm:generalconv}, we take the approach in \citet{ouyang2015accelerated}, and analyze GeNI-ADMM by viewing \cref{eq:ADMMProb} through its formulation as saddle point problem \cref{eq:saddlepoint}.  
Let $W = X \times Z \times U$, $\hat{w} = (\hat{x},\hat{z},\hat{u})$, and $w = (x,z,u)$, where $w, \hat w \in W$. 
Define the \emph{gap function}
\begin{equation}
\label{eq:GapFun}
	Q(\hat w, w) \coloneqq [f(x)+g(z)+\langle\rho\hat{u},Mx+Nz\rangle]-[f(\hat{x})+g(\hat{z})+\langle \rho u, M\hat{x}+N\hat z \rangle].
\end{equation}
The gap function naturally arises from the saddle-point formulation \cref{eq:saddlepoint}. 
By construction, it  
is concave in its first argument and convex in the second, 
and satisfies the important inequality
\[
Q(w^\star,w) = L_\rho(x,z,u^\star)-L_\rho(x^\star,z^\star,u) \geq 0, 
\]
which follows by definition of $(x^\star,z^\star,u^\star)$ being a saddle-point. 
Hence the gap function may be viewed as measuring the distance to the saddle $w^\star$.

Further, given a closed set $U \subset \mathcal Z$ and $v\in \mathcal{Z}$ we define
\begin{align}
\label{eq:lU-def}
\ell_U(v, w) &\coloneqq \sup_{\hat{u} \in U} 
\{Q(\hat w, w) + \langle v, \rho \hat{u}\rangle \mid
\hat{w} = (\hat{x},\hat{z},\hat{u}),~\hat x = x^\star,~\hat z = z^\star\}\\
&= f(x)+g(z)-p^\star +\sup_{\hat u\in U}\langle \rho \hat u, v-(Mx+Nz) \rangle. \nonumber 
\end{align}
The following lemma of \citet{ouyang2015accelerated} relates $Q(\hat w, w)$ to the suboptimality and feasibility gaps.
\begin{lemma}
\label{lemma:GapBnd}
For any $U \subseteq \mathcal Z$, suppose $\ell_U(Mx + Nz, w) \le \epsilon < \infty$ and $\|Mx +Nz\| \le \delta$, where $w = (x, z, u) \in W$.
Then
\begin{equation}
    f(x)+g(z)-p^{\star}\leq \epsilon. 
\end{equation}
In other words, 
$(x, z)$ is an approximate solution of \eqref{eq:ADMMProb} with suboptimality gap $\epsilon$ and feasibility gap~$\delta$. 
Further, if $U = \mathcal Z$, for any $v$ such that $\ell_U(v, w) \le \epsilon < \infty$ and $\|v\| \le \delta$, we have $v = Mx - z$. 
\end{lemma}
\cref{lemma:GapBnd} shows that if we can find $w$ such that $\ell_U(Mx + Nz, w)\leq \epsilon$ and $\|Mx + Nz\| \le \delta$, 
then we have an approximate optimal solution to \eqref{eq:ADMMProb} with gaps $\epsilon$ and $\delta$, that is, $\ell_{U}(Mx+Nz,w)$ controls the suboptimality and feasibility gaps.

\subsection{Controlling the gap function}
\cref{lemma:GapBnd} shows the key to establishing convergence of GeNI-ADMM is to achieve appropriate control over the gap function.
To accomplish this, we use the optimality conditions of the $x$ and $ z$ subproblems. 
However, as the subproblems are only solved approximately, the inexact solutions satisfy \emph{perturbed} optimality conditions.
To be able to reason about the optimality conditions under inexact solutions, the iterates must remain bounded. 
Indeed, if the iterates are unbounded, they can fail to satisfy the subproblem optimality conditions arbitrarily badly.
Fortunately, the following proposition shows the iterates remain bounded.

\begin{proposition}[GeNI-ADMM iterates remain bound]
\label{prop:bnded_iterates}
    Let $\{\varepsilon_x^k\}_{k\geq 1}, \{\varepsilon_z^k\}_{k\geq 1}, \{\Theta^k\}_{k\geq 1}, \{\Psi^k\}_{k\geq 1}$, and $\rho>0$ be given. 
    Instate Assumptions \ref{assp:existence}-\ref{assp:metric_evol}. 
    Run \cref{alg-gen_newton_admm}, then the output sequences $\{\tilde x^k\}_{k\geq 1}, \{\tilde z^k\}_{k\geq 1}, \{\tilde u^k\}_{k\geq 1}$ are bounded.
    That is, there exists $R>0$, such that
    \[
    \{\tilde x^k\}_{k\geq 1}\subset B(x^\star,R), \quad \{\tilde z^k\}_{k\geq 1}\subset B(z^\star,R), \quad \{\tilde u^k\}_{k\geq 1}\subset B(u^\star,R). 
    \]
\end{proposition}
The proof of \cref{prop:bnded_iterates} is provided in \cref{subsection:SubLinConvPfs}.

As the iterates remain bounded, we can show that the optimality conditions are approximately satisfied at each iteration.  
The precise form of these perturbed optimality conditions is given in \cref{lemma:inexactxopt,lemma:inexactzopt}.
Detailed proofs establishing these lemmas are given in \cref{subsection:SubLinConvPfs}. 

 \begin{lemma}[Inexact $x$-optimality condition]
\label{lemma:inexactxopt}
Instate Assumptions \ref{assp:existence}-\ref{assp:metric_evol}.
Suppose $\tilde x^{k+1}$ is an $\varepsilon^k_x$-minimizer of the $x$-subproblem under Assumption \ref{assp:x-subprob-orc}. 
Then for some absolute constant $C_x>0$ we have 
\begin{equation}
    \label{eq:inexactxoptsimplified}
    \begin{aligned}
    &\langle \nabla f_1(\tilde x^{k+1})+\nabla f_2(\tilde x^k), \tilde x^{k+1}-x^\star\rangle \leq \langle \Theta^k(\tilde x^{k+1}-\tilde x^k),x^\star-\tilde x^{k+1}\rangle\\
    &+\rho \langle M\tilde x^{k+1}+N\tilde z^k+\tilde u^k,M(x^\star-\tilde x^{k+1})\rangle+C_x\varepsilon_x^k.
    \end{aligned}
\end{equation}
\end{lemma}
\begin{lemma}[Inexact $z$-optimality condition]
\label{lemma:inexactzopt}
Instate Assumptions \ref{assp:existence}-\ref{assp:metric_evol}.
Suppose $\tilde z^{k+1}$ is an $\varepsilon^k_z$-minimum of the $z$-subproblem under \cref{def:TypeInexactness}. 
Then for some absolute constant $C_z>0$ we have 
\begin{equation}
\label{eq:inexactzoptsimplified}
 g(z^\star)-g(\tilde z^{k+1})\geq \langle -\rho N^{T}\tilde u^{k+1}+\Psi^k(\tilde z^k-\tilde z^{k+1}), z^{\star}-\tilde z^{k+1}\rangle-C_z\sqrt{\varepsilon_z^k}.
\end{equation}
\end{lemma}

When $\varepsilon_x^k = \varepsilon_z^k = 0$ the approximate optimality conditions of \cref{lemma:inexactxopt,lemma:inexactzopt} collapse to the exact optimality conditions. 
We also note while \cref{lemma:inexactxopt} (\cref{lemma:inexactzopt}) necessarily holds when $\tilde x^{k+1}$ ($\tilde z^{k+1})$ is an $\varepsilon_x^{k}$-approximate minimizer ($\varepsilon_z^k$-minimum), the converse does not hold.
With \cref{lemma:inexactxopt,lemma:inexactzopt} in hand, we can establish control of the gap function for one iteration. 

\begin{lemma}
\label{lemma:GapBndInexactCase}
Instate Assumptions \ref{assp:existence}-\ref{assp:metric_evol}..
Let $\tilde w^{k+1} = (\tilde x^{k+1}, \tilde z^{k+1}, \tilde u^{k+1})$ denote the iterates generated by \cref{alg-gen_newton_admm} at iteration $k$. Set $w = (x^\star,z^\star,u)$, then
the gap function $Q$ satisfies
\begin{align*}
   	Q(w,\tilde w^{k+1}) & \leq \frac{1}{2}\left(\|\tilde x^k-x^\star\|^2_{\Theta^k}-\|\tilde x^{k+1}-x^\star\|^2_{\Theta^k}\right) + 
	\frac{1}{2} \left(\|\tilde z^k - z^\star\|_{\Psi^k+\rho N^TN}^2 - \|\tilde z^{k+1} - z^\star\|_{\Psi^k+\rho N^TN}^2\right)\\
    &+\frac{\rho}{2}\left(\|\tilde u^k-u\|^2-\|\tilde u^{k+1}-u\|^2\right)+C_x\varepsilon_x^k+C_z\sqrt{\varepsilon_z^k}. 
\end{align*}
\end{lemma}
\begin{proof}
From the definition of $Q$,
\begin{align*}
    Q(w,\tilde{w}^{k+1}) = f(\tilde x^{k+1})-f(x^\star)+g(\tilde z^{k+1})-g(z^\star)+\langle \rho u, M\tilde{x}^{k+1}+N\tilde{z}^{k+1}\rangle -\langle \rho \tilde u^{k+1},Mx^\star+Nz^\star\rangle.
\end{align*}
Our goal is to upper bound $Q(w,\tilde{w}^{k+1})$. 
We start by bounding $f(\tilde x^{k+1})-f(x^\star)$ as follows:
\begin{align*}
    &f(\tilde x^{k+1})-f(x^\star) = \left(f_1(\tilde x^{k+1})-f_1(x^\star)\right)+\left(f_2(\tilde x^{k+1})-f_2(\tilde x^k)\right)+\left(f_2(\tilde x^k)-f_2(x^\star)\right)\\
    &\overset{(1)}{\leq} \langle \nabla f_1(\tilde x^{k+1}),\tilde x^{k+1}-\tilde x^\star\rangle+ \langle \nabla f_1(\tilde x^{k}),\tilde x^{k+1}-\tilde x^k\rangle +\frac{1}{2}\|\tilde x^{k+1}-\tilde x^k\|^2_{\Theta_k}+f_2(\tilde x^k)-f_2(x^\star)\\
    &\overset{(2)}{\leq} \langle \nabla f_1(\tilde x^{k+1}),\tilde x^{k+1}-\tilde x^\star\rangle+\langle \nabla f_2(\tilde x^k),\tilde x^{k+1}-\tilde x^k\rangle +\frac{1}{2}\|\tilde x^{k+1}-\tilde x^k\|^2_{\Theta_k}+\langle \nabla f_2(\tilde x^k),\tilde x^k-x^\star\rangle \\
    &= \langle \nabla f_1(\tilde x^{k+1})+\nabla f_2(\tilde x^k),\tilde x^{k+1}-x^\star\rangle+\frac{1}{2}\|\tilde x^{k+1}-\tilde x^{k}\|^2_{\Theta^k},
\end{align*}
where $(1)$ uses convexity of $f_1$ and $1$-relative smoothness of $f_2$, and $(2)$ uses convexity of $f_2$.
Inserting the upper bound on $f(\tilde x^{k+1})-f(x)$ into the expression for $Q(w,\tilde{w}^{k+1})$, we find
\begin{equation}
\label{eq:GapBnd2}
\begin{aligned}
    Q(w,\tilde w^{k+1})\leq& \langle \nabla f_1(\tilde x^{k+1})+\nabla f_2(\tilde x^k), \tilde x^{k+1}-x^\star\rangle +\frac{1}{2}\|\tilde x^{k+1}-\tilde x^k\|^2_{\Theta^k}+g(\tilde z^{k+1})-g(z^\star)\\
    &+\langle \rho u, M\tilde x^{k+1}-N\tilde z^{k+1}\rangle-\langle \rho \tilde u^{k+1}, Mx^\star-Nz^\star\rangle. 
\end{aligned}
\end{equation}
Now, using the inexact optimality condition for the $x$-subproblem (\cref{lemma:inexactxopt}), the above display becomes
\begin{align}
    & Q(w,\tilde w^{k+1})\leq \langle \Theta^k(\tilde x^{k+1}-\tilde x^k),x^\star-\tilde x^{k+1}\rangle+\rho \langle M\tilde x^{k+1}+N\tilde z^k+\tilde u^k,M(x^\star-\tilde x^{k+1})\rangle +C_x\varepsilon_x^k\\ 
    &+\frac{1}{2}\|\tilde x^{k+1}-\tilde x^k\|^2_{\Theta^k}+g(\tilde z^{k+1})-g(z^\star)+\langle \rho u, M\tilde x^{k+1}+N\tilde z^{k+1}\rangle-\langle \rho \tilde u^{k+1}, Mx^\star+Nz^\star\rangle. 
\end{align}
Similarly, applying the inexact optimality for the $z$-subproblem (\cref{lemma:inexactzopt}), we further obtain
\begin{align}
\label{eq:int_gap_bnd}
    & Q(w,\tilde w^{k+1})\leq \langle \Theta^k(\tilde x^{k+1}-\tilde x^k),x^\star-\tilde x^{k+1}\rangle+\rho \langle M\tilde x^{k+1}+N\tilde z^k+\tilde u^k,M(x^\star-\tilde x^{k+1})\rangle +C_x\varepsilon_x^k\nonumber \\ 
    &+\frac{1}{2}\|\tilde x^{k+1}-\tilde x^k\|^2_{\Theta^k}-\rho \langle N^{T}\tilde u^{k+1},\tilde z^{k+1}-z^\star\rangle +\langle \Psi^k(\tilde z^{k+1}-\tilde z^k),z^\star-\tilde z^{k+1}\rangle+C_z\sqrt{\varepsilon_z^k} \\
    &+\langle \rho u, M\tilde x^{k+1}+N\tilde z^{k+1}\rangle-\langle \rho \tilde u^{k+1}, Mx^\star+Nz^\star\rangle \nonumber. 
\end{align}
We now simplify \eqref{eq:int_gap_bnd} by combining terms.
Some basic manipulations show the terms on line 2 of \eqref{eq:int_gap_bnd} may be rewritten as
\begin{align*}
    &-\rho \langle N^{T}\tilde u^{k+1},\tilde z^{k+1}-z^\star\rangle +\rho \langle u, M\tilde x^{k+1}+N\tilde z^{k+1}\rangle-\rho \langle \tilde u^{k+1},Mx^\star+Nz^\star\rangle \\
    &= \rho \langle u-\tilde u^{k+1},\tilde M\tilde x^{k+1}+N\tilde z^{k+1}\rangle
    -\rho \langle \tilde u^{k+1},M(x^\star-\tilde x^{k+1}) \rangle.
\end{align*}
We can combine the preceding display with the second term of line 1 in \eqref{eq:int_gap_bnd} to reach
\begin{align*}
  &\rho \langle M\tilde x^{k+1}+N\tilde z^k+\tilde u^k,M(x^\star-\tilde x^{k+1})\rangle
  -\rho \langle \tilde u^{k+1},M(x^\star-\tilde x^{k+1}) \rangle +\rho \langle u-\tilde u^{k+1},M\tilde x^{k+1}+N\tilde z^{k+1}\rangle \\
  &= \rho\langle \tilde u^{k+1}+N(\tilde z^{k}-\tilde z^{k+1}),M(x^\star-\tilde x^{k+1})\rangle
  -\rho \langle \tilde u^{k+1},M(x^\star-\tilde x^{k+1}) \rangle+ \rho \langle u-\tilde u^{k+1},M\tilde x^{k+1}+N\tilde z^{k+1}\rangle\\
  &= \rho \langle N(\tilde z^{k}-\tilde z^{k+1}),M(x^\star-\tilde x^{k+1})+ \rho \langle u-\tilde u^{k+1},M\tilde x^{k+1}+N\tilde z^{k+1}\rangle \\
  & = \rho \langle \Nm(\tilde z^{k+1}-\tilde z^{k}),M(x^\star-\tilde x^{k+1})\rangle +\rho \langle \tilde u^{k+1}-u,\tilde u^k-\tilde u^{k+1}\rangle, 
\end{align*}
where $\Nm = -N.$
Inserting the preceding simplification into \eqref{eq:int_gap_bnd}, we reach 
\begin{equation}
\label{eq:simpli_gap_bnd}
\begin{aligned}
    Q(w,\tilde w^{k+1}) & \le
    \langle \Theta^k(\tilde x^{k+1}-\tilde x^k),x^\star-\tilde x^{k+1}\rangle + \rho \langle \tilde u^{k+1}-u, \tilde u^{k}-\tilde u^{k+1}\rangle \\
    &+ \langle \Psi^k(\tilde z^{k+1}-\tilde z^k),z^\star-\tilde z^{k+1}\rangle + \rho \langle \Nm(\tilde z^{k+1}-\tilde z^{k}),M(x^\star-\tilde x^{k+1})\rangle\\
    &+\frac{1}{2} \|\tilde x^{k+1} - \tilde x^k\|_{\Theta^k}^2+C_x\varepsilon_x^k+C_z\sqrt{\varepsilon^k_z}.
\end{aligned}
\end{equation}
Now, we bound the first two leading terms in line 1 of \eqref{eq:simpli_gap_bnd} by invoking the identity $\langle a-b,\Upsilon(c-d)\rangle = 1/2\left(\|a-d\|^2_\Upsilon-\|a-c\|^2_\Upsilon\right)+1/2\left(\|c-b\|^2_\Upsilon-\|d-b\|^2_\Upsilon\right)$ to obtain
\begin{align*}
    &\langle \Theta^k(\tilde x^{k+1}-\tilde x^k),x^\star-\tilde x^{k+1}\rangle+\rho\langle \tilde u^{k+1}-u, \tilde u^k-\tilde u^{k+1}\rangle = \\
    &\frac{1}{2}\left(\|x^\star-\tilde x^k\|^2_{\Theta^k}-\|x^\star-\tilde x^{k+1}\|^2_{\Theta^k}\right)-\frac{1}{2}\|\tilde x^{k+1}-\tilde x^k\|^2_{\Theta^k} + \frac{\rho}{2}\left(\|\tilde u^k-u\|^2-\|\tilde u^{k+1}-u\|^2-\|\tilde u^k-\tilde u^{k+1}\|^2\right).
\end{align*}
Similarly, to bound the third and fourth terms in \eqref{eq:simpli_gap_bnd}, we again invoke $(a-b)^T\Upsilon(c-d) = 1/2\left(\|a-d\|^2_\Upsilon-\|a-c\|^2_\Upsilon\right)+1/2\left(\|c-b\|^2_\Upsilon-\|d-b\|^2_\Upsilon\right)$ which yields
\begin{align*}
    \langle \Psi^k(\tilde z^{k+1}-\tilde z^k),z^\star-\tilde z^{k+1}\rangle = \frac{1}{2}\left(\|\tilde z^k-z^\star\|^2_{\Psi^k}-\|\tilde z^{k+1}-z^\star\|^2_{\Psi^k}\right)-\frac{1}{2}\|\tilde z^{k+1}-\tilde z^k\|_{\Psi^k}^2,
\end{align*}
\begin{align*}
&\rho\langle \Nm(\tilde z^k-\tilde z^{k+1}),M(\tilde x^{k+1}-x^\star)\rangle\\
& = \frac{\rho}{2} \left(\|\Nm \tilde z^k - M x^\star \|^2 - \|\tilde z^{k+1} - M x^\star \|^2 + \|\Nm\tilde z^{k+1} - M \tilde x^{k+1}\|^2 - \|\Nm\tilde z^k - M \tilde x^{k+1}\|^2\right) \\
& = \frac{\rho}{2} \left(\|N\tilde z^k-Nz^\star \|^2 - \|N\tilde z^{k+1}-Nz^\star\|^2  - \|M \tilde x^{k+1}+N\tilde z^k\|^2\right) + \frac{\rho}{2} \|\tilde u^k - \tilde u^{k+1}\|^2.
\end{align*}
Putting everything together, we conclude
\begin{align*}
	Q(w,\tilde w^{k+1}) & \le \frac{1}{2}\left(\|\tilde x^k-x^\star\|^2_{\Theta^k}-\|\tilde x^{k+1}-x^\star\|^2_{\Theta^k}\right) + 
	\frac{1}{2} \left(\|\tilde z^k - z^\star\|_{\Psi^k+\rho N^TN}^2 - \|\tilde z^{k+1} - z^\star\|_{\Psi^k+\rho N^TN}^2\right)\\
    &+\frac{\rho}{2}\left(\|\tilde u^k-u\|^2-\|\tilde u^{k+1}-u\|^2\right)+C_x\varepsilon_x^k+C_z\sqrt{\varepsilon_z^k}. 
\end{align*}
as desired.
\end{proof}

\subsection{Proof of \cref{thm:generalconv}}
With \cref{lemma:GapBndInexactCase} in hand, we are ready to prove \cref{thm:generalconv}. 
\proof
From \Cref{lemma:GapBndInexactCase}, we have for each $k$ that
\begin{align*}
	Q(w,\tilde w^{k+1}) & \le \frac{1}{2}\left(\|\tilde x^k-x^\star\|^2_{\Theta^k}-\|\tilde x^{k+1}-x^\star\|^2_{\Theta^k}\right) + 
	\frac{1}{2} \left(\|\tilde z^k - z^\star\|_{\Psi^k+\rho N^TN}^2 - \|\tilde z^{k+1} - z^\star\|_{\Psi^k+\rho N^TN}^2\right)\\
    &+\frac{\rho}{2}\left(\|\tilde u^k-u\|^2-\|\tilde u^{k+1}-u\|^2\right)+C_x\varepsilon_x^k+C_z\sqrt{\varepsilon_z^k}. 
\end{align*}
Now, summing up the preceding display from $k=1$ to $t$ and using $w = (x^\star,z^\star, u)$, we obtain
\begin{align*}
    \sum^{t+1}_{k=2}Q(x^\star,z^\star, u;\tilde{w}^{k}) & \leq \frac{1}{2}\underbrace{\sum^{t}_{k=1} \left(\|\tilde x^k-x^\star\|^2_{\Theta^k}-\|\tilde x^{k+1}-x^\star\|^2_{\Theta^k}\right)}_{T_1}\\
    &+\underbrace{\frac{\rho}{2}\sum^{t}_{k=1}\left(\|\tilde z^k - z^\star\|_{\Psi^k+\rho N^TN}^2 - \|\tilde z^{k+1} - z^\star\|_{\Psi^k+\rho N^TN}^2\right)}_{T_2} \\
    & + \underbrace{\frac{\rho}{2}\sum^{t}_{k=1} \left(\|\tilde u^k-u\|^2-\|\tilde u^{k+1}-u\|^2\right)}_{T_3}+C_x\mathcal E_x+C_z\mathcal E_z.
\end{align*}

We now turn to bounding $T_1$. Using the definition of $T_1$, we find
\begin{align*}
    T_1 & = \frac{1}{2}\sum^{t}_{k=1} \left(\|\tilde x^k-x^\star\|^2_{\Theta^k}-\|\tilde x^{k+1}-x^\star\|^2_{\Theta^k}\right) \\ 
    & = \frac{1}{2}\left(\|\tilde x^1-x^\star\|^2_{\Theta^1}-\|\tilde x^{t+1}-x^\star\|^2_{\Theta^t}\right)
    +\frac{1}{2}\sum_{k=2}^{t}\left(\|\tilde x^k-x^\star\|^2_{\Theta^k}-\|\tilde x^k-x^\star\|^2_{\Theta^{k-1}}\right).
\end{align*}
Now using our hypotheses on the sequence $\{\Theta^k\}_k$ in \eqref{eq:Hcond}, we obtain
\begin{align*}
   \|\tilde x^k -x^\star\|^2_{\Theta^k}-\|\tilde x^k-x^\star\|^2_{\Theta^{k-1}} &= (\tilde x^k -x^\star)^{T}(\Theta^k-\Theta^{k-1})(\tilde x^k -x^\star)\leq \zeta^{k-1} \| \tilde x^k-x^\star\|^2_{\Theta^{k-1}}\\
   &\leq \tau_\zeta \zeta^{k-1}\| \tilde x^k-x^\star\|^2_{\Theta^{1}} .
\end{align*}
Inserting the previous bound into $T_1$, we reach

\begin{align*}
    T_1 &\leq \frac{1}{2}\|\tilde x^1-x^\star\|^2_{\Theta^1}+\frac{1}{2}\tau_\zeta \sum_{k=1}^{t}\zeta^{k-1}\|\tilde x^{k}-x^{\star}\|^2_{\Theta^1} - \frac{1}{2}\|\tilde x^{t+1}-x^\star\|^2_{\Theta^t} \\ 
    &\leq \frac{1}{2}\left(\|\tilde x^1-x^\star\|^2_{\Theta^1}+C_\zeta\mathcal{E}_{\zeta} \right) - \frac{1}{2}\|\tilde x^{t+1}-x^\star\|^2_{\Theta^t}\\
    & \leq \frac{1}{2}\left(d^2_{x^\star, \Theta^1}+C_\zeta \mathcal{E}_{\zeta} \right).
\end{align*}

Next, we bound $T_2$. 
\begin{align*}
    T_2 & = \frac{1}{2}\sum^{t}_{k=1}\left(\|\tilde z^k - z^\star\|^2 - \|\tilde z^{k+1} - z^\star\|^2 \right) + \frac{1}{2}\sum_{k=1}^{t} \left(\|\tilde z^k - z^\star\|_{\Psi^k+\rho N^TN}^2 - \|\tilde z^{k+1} - z^\star\|_{\Psi^k+\rho N^TN}^2\right) \\
    & \leq \frac{\rho}{2}\sum^{t}_{i=1} \left(\|\tilde z^k-z^\star \|^2 - \|\tilde z^{k+1} - z^\star \|^2\right) + \frac{1}{2}\left(\|\tilde z^1-z^\star\|^2_{\Psi^1}+C_\zeta \mathcal E_\zeta \right)\\
    &= \frac{1}{2}\|\tilde z^1-z^\star\|_{\Psi^1+\rho I}^2 = \frac{1}{2}d_{z^\star,\Psi^1_{\rho,N}}^{2}+\frac{C_\zeta}{2} \mathcal E_\zeta.
\end{align*}

Last, $T_3$ is a telescoping sum, hence
\begin{align*}
    T_2 & = \frac{\rho}{2}\sum^{t}_{k=1} \left(\|\tilde u^k-u\|^2-\|\tilde u^{k+1}-u\|^2\right) = \frac{\rho}{2}\left(\|\tilde u^1-u\|^2-\|\tilde u^{t+1}-u\|^2\right). 
\end{align*}

Using our bounds on $T_1$ through $T_3$, we find
\begin{equation*}
\begin{aligned}
    \sum^{t+1}_{k=2}Q(x^\star,z^\star, u;\tilde w^{k}) & \leq \frac{1}{2}\left(d^2_{x^\star,\Theta^1}+C_\zeta\mathcal{E}_{\zeta}\right)+\frac{1}{2}d_{z^\star,\Psi^1_{\rho,N}}^{2}+\frac{1}{2}C_\zeta \mathcal E_\zeta+ C_x\mathcal{E}_x+C_z\mathcal{E}_z\\ & + \frac{\rho}{2}\left(\|\tilde u^1-u\|^2-\|\tilde u^{t+1}-u\|^2\right) \\
    &= \frac{1}{2}d^2_{x^\star,\Theta^1}+\frac{1}{2}d_{z^\star,\Psi^1_{\rho,N}}^{2}+ C_x\mathcal{E}_x+C_z\mathcal{E}_z+C_\zeta\mathcal{E}_{\zeta}+ \frac{\rho}{2}\left(\|\tilde u^1-u\|^2-\|\tilde u^{t+1}-u\|^2\right),\\
\end{aligned}
\end{equation*}
Now, as $\bar w^{t+1} = \frac{1}{t}\sum^{t+1}_{k=2}\tilde w^{k}$, the convexity of $Q$ in its second argument yields
\begin{equation}
\label{eq:qconvexity}
\begin{aligned}
    Q(x^\star,z^\star, u;\bar w^{t+1}) & \leq \frac{1}{t}\sum_{k=2}^{t+1}Q(x^\star,z^\star, u;\tilde w^{k})\\
    &\leq \frac{1}{t}\left(\frac{1}{2}d^2_{x^\star,\Theta^1}+\frac{\rho}{2}d_{z^\star,\Psi_{\rho,N}^1}^{2}+ C_x\mathcal{E}_x+C_z\mathcal{E}_z+C_\zeta\mathcal{E}_{\zeta}+ \frac{\rho}{2}\left(\|\tilde u^1-u\|^2-\|\tilde u^{t+1}-u\|^2\right)\right).
\end{aligned}
\end{equation}
Define $\Gamma \coloneqq \frac{1}{2}d^2_{x^\star,\Theta^1}+\frac{1}{2}d^2_{z^\star,\Psi^1_{\rho,N}}+ C_x\mathcal{E}_x+C_z\mathcal{E}_z+C_\zeta\mathcal{E}_{\zeta}$.
Since $Q(w^\star, \bar w^{t+1}) \ge 0$, by \eqref{eq:qconvexity} we reach 
\begin{equation*}
    \|\tilde u^{t+1}-u^\star\|^2 \le \frac{2}{\rho}\Gamma+d^2_{u^\star}.
\end{equation*}
Let $\tilde v^{t+1} = \frac{1}{t}\left(\tilde u^1 - \tilde u^{t+1}\right)$. 
Then we can bound $\|\tilde v_{t+1}\|^2$ as 
\begin{equation*}
\begin{aligned}
 \|\tilde v^{t+1}\|^2 \le & \frac{2}{t^2}\left(\|\tilde u^1-u^\star\|^2+\|\tilde u^{t+1}-u^\star\|^2\right) \leq \frac{4}{t^2}\left(\frac{\Gamma}{\rho}+d^2_{u^\star}\right).
\end{aligned}
\end{equation*}
By \eqref{eq:qconvexity}, given the fact $\tilde u^1 = 0$, we also have 
\begin{equation*}
\begin{aligned}
    Q(x^\star,z^\star, u;\bar w^{t+1}) 
    &\leq \frac{\Gamma-\rho \langle \tilde u^1 - \tilde u^{t+1}, u\rangle}{t} = \frac{\Gamma}{t}-\rho \langle \tilde v^{t+1}, u\rangle,
\end{aligned}
\end{equation*}
where the equality follows from the definition of $\tilde v^{t+1}$. 
Hence for any $u$
\begin{align*}
Q(x^\star,z^\star, u;\bar w^{t+1})+\langle \tilde v^{t+1},\rho u \rangle 
&\leq  \frac{\Gamma}{t},
\end{align*}
and therefore 
\begin{equation*}
    \ell_U(\tilde v^{t+1}, \bar w^{t+1}) \leq \frac{1}{t}\left(\frac{1}{2}d^2_{x^\star,\Theta^1}+\frac{1}{2}d_{z^\star,\Psi^1_{\rho,N}}^{2}+ C_x\mathcal{E}_x+C_z\mathcal{E}_z+C_\zeta\mathcal{E}_{\zeta}\right).
\end{equation*}
We finish the proof by invoking \cref{lemma:GapBnd}. $\Box$ 

\section{Linear convergence of GeNI-ADMM}
\label{section:linearconvergence}
In this section, we seek to establish linear convergence results 
for \cref{alg-gen_newton_admm}. 
In general, the linear convergence of ADMM relies on strong convexity of the objective function \citep{MAL-016,nishihara2015general,parikh2014proximal}. Consistently, the linear convergence of GeNI-ADMM also requires strong convexity. Many applications of GeNI-ADMM fit into this setting, such as elastic net \citep{friedman2010regularization}. However, linear convergence is not restricted to strongly convex problems. It has been shown that local linear convergence of ADMM can be guaranteed even without strong convexity \citep{yuan2020discerning}. Experiments in \cref{section:experiments} show the same phenomenon for GeNI-ADMM: it converges linearly after a couple of iterations when the iterates reach some manifold containing the solution. 
The linear convergence theory of GeNI-ADMM provides a way to understand this phenomenon. 
We first list the additional assumptions required for linear convergence:
\begin{assumption}[Optimization is over the whole space]
\label{assp:opt_space}
    The sets $X$ and $Z$ in \eqref{eq:ADMMProb} satisfy 
    \[
    X = \mathcal X,~ \text{and}~ Z = \mathcal Z.
    \]
\end{assumption}
Assumption \ref{assp:opt_space} states that the optimization problem in \eqref{eq:ADMMProb} is over the entire spaces $\mathcal X$ and $\mathcal Z$, not closed subsets.
This assumption is met in many practical optimization problems of interest, 
where $\mathcal X = \mathcal Z = \mathcal H = \mathbb{R}^{d}$.
Moreover, it is consistent with prior analyses such as \citet{deng2016global}, who specialize their analysis to the setting of the last sentence. 

\begin{assumption}[Regularity of $f$]
\label{assp:linearregularity}
    The function $f$ is finite valued, strongly convex with parameter $\sigma_f$, and smooth with parameter $L_f$. 
\end{assumption}
Assumption \ref{assp:linearregularity} imposes standard regularity conditions on $f$, in addition to the conditions of Assumption \ref{assp:regularity}. 
\begin{assumption}[Non-degeneracy of constraint operators]
\label{assp:fullrank}
    The linear operators $MM^T$ and $N^TN$ are invertible.
\end{assumption}
Assumption \ref{assp:fullrank} is consistent with prior analyses of ADMM-type schemes under strong convexity, 
such as \citet{deng2016global}, who make this assumption in their analysis of Generalized ADMM. 
Moreover, the assumption holds in many important application problems, 
especially those that arise from machine learning 
where, typically, $M=I$ and $N=-I$. 
\begin{assumption}[Geometric decay of the forcing sequences]
\label{assp:geodecay}
   There exists a constant $q>0$ such that the forcing sequences $\{\varepsilon^k_x\}_{k \ge 1}$, $\{\varepsilon^k_z\}_{k \ge 1}$ satisfy 
   \begin{equation}
       \label{eq:geodecay}
       \varepsilon^{k+1}_x \le \varepsilon^k_x / (1 + q), \; \text{and} \;\;\varepsilon^{k+1}_z \le \varepsilon^k_z / (1 + q)^2.
   \end{equation}
   Moreover, we assume \cref{alg-gen_newton_admm} is equipped with oracles for solving the $x$ and $z$-subproblems, which at each iteration produce approximate solutions $\tilde{x}^{k+1}$, $\tilde z^{k+1}$ satisfying:
\begin{align*}
    & \tilde x^{k+1} \xepsappx  \underset{x \in X}{\argmin}\{ f_1(x) + \langle \nabla f_2(\tilde x^k ), x - \tilde x^k\rangle+\frac{1}{2}\|x - \tilde x^k\|^2_ {\Theta^k} + \frac{\rho}{2}\|Mx+N\tilde z^k + \tilde u^k\|^2\},\\
    & \tilde z^{k+1} \zfepsappx \underset{z\in Z}{\argmin} \{g(z)+ \frac{\rho}{2}\|M\tilde x^{k+1} +Nz + \tilde u^k\|^2+\frac{1}{2}\|z-\tilde z^k\|^2_{\Psi^k}\}. \end{align*}
\end{assumption}

Assumption \ref{assp:geodecay} replaces Assumption \ref{assp:x-subprob-orc} and requires the inexactness sequences to decay geometrically. 
Compared with the sublinear convergence result, linear convergence requires more accurate solutions to the subproblems. 
Again, since the $z$-subproblem inexactness is weaker than the $x$-subproblem inexactness, $\{\varepsilon^{k}_z\}_{k \ge 1}$ should have a faster decay rate $(1+q)^2$ than the decay rate ($1+q$) of $\{\varepsilon^{k}_x\}_{k \ge 1}$.

The requirement that the forcing sequences decay geometrically is somewhat burdensome, as it leads to the subproblems needing to be solved to higher accuracy sooner than if the forcing sequences were only summable.
\ifpreprint
Fortunately, this condition seems to be an artifact of the analysis;  our numerical experiments with strongly convex $f$ (\cref{subsection:lin_exp}) only use summable forcing sequences but show linear convergence of GeNI-ADMM.
\fi

\subsection{Our approach}
Inspired by \citet{deng2016global}, we take a $\emph{Lyapunov function}$ approach to proving linear convergence.
Let $\tilde w = (\tilde x,\tilde z, \tilde u)$, and $w^\star = (x^\star, z^\star, u^\star)$. 
We define the Lyapunov function: 
\[
\Phi^k = \frac{1}{\rho}\|\tilde x^k-x^\star\|_{\Theta^k}^2+\frac{1}{\rho}\|\tilde z^k-\tilde z^\star\|_{\Psi^k+\rho N^TN}^2+\|\tilde u^k-u^\star\|^2 = \|\tilde w^k-w^\star\|_{G^k}^2, 
\]
where
\[G^k \coloneqq \begin{pmatrix}
    \frac{1}{\rho}\Theta^k & 0 & 0 \\
    0 & \frac{1}{\rho} \Psi^k+N^TN & 0 \\
    0 & 0 & I
\end{pmatrix}.
\]
Our main result in this section is the following theorem, which shows the Lyapunov function converges linearly to $0$.
\begin{theorem}\label{thm:linearconvabs}
Instate Assumptions \ref{assp:existence}-\ref{assp:regularity}, and Assumptions \ref{assp:metric_evol}-\ref{assp:geodecay}. 
Moreover, suppose that $\theMin$ in \eqref{eq:mat_reg_assmp} satisfies $\theMin>L^2_f/\sigma_f$.
Then there exist constants $\delta$ and $S > 0$ such that if 
$q > \delta$, 
	\begin{equation}\label{eq:linearconvabs}
		(1+\delta)^{k}\Phi^{k} \le \tau_\zeta \Phi^1 + S.
	\end{equation}
    Hence after $k=\mathcal O\left(\frac{1}{\delta}\log\left(\frac{\tau_\zeta \Phi^1}{\epsilon}\right)\right)$ iterations, 
    \[
    \Phi^k\leq \epsilon. 
    \]
\end{theorem}
The proof of \cref{thm:linearconvabs} is deferred to \cref{subsection:LinConPf}.
As $G^k\succ 0$ for all $k$, \cref{thm:linearconvabs} implies the iterates $(\tilde x^k, \tilde z^k,\tilde u^k)$ converge linearly to optimum $(x^\star,z^\star,u^\star)$.
Thus, despite inexactly solving approximations of the original ADMM subproblems, GeNI-ADMM still enjoys linear convergence when the objective is strongly convex. 
An unattractive aspect of \cref{thm:linearconvabs} is the requirement that $\theMin$ satisfy $\theMin>L_f^2/\sigma_f$.
Although this condition can be enforced by adding a damping term $\sigma I$ to $\Theta^k$, this is undesirable as a large regularization term can lead GeNI-ADMM to converge slower.  
We have found that this condition is unnecessary---our numerical experiments do not enforce this condition, yet GeNI-ADMM achieves linear convergence.
We believe this condition to be an artifact of the analysis, which stems from lower bounding the term $-\langle \nabla f_2(\tilde x^{k+1})-\nabla f_2(\tilde x^k),\tilde x^{k+1}-x^\star\rangle$ in \cref{lemma:suffdescent}.

\subsection{Sufficient descent}
From \cref{thm:linearconvabs}, to establish linear convergence of GeNI-ADMM, it suffices to show $\Phi^k$ decreases geometrically. 
We take two steps to achieve this. 
First, we show that $\|\tilde w^k-w^\star\|_{G^k}^2-\|\tilde w^{k+1}-w^\star\|_{G^k}^2$ decreases for every iteration $k$ (\cref{lemma:suffdescent}). 
Second, we show that $\Phi^{k+1}$ decreases geometrically by a factor of $1/(1+\delta)$ with respect to $\Phi^k$ and some small error terms that stem from inexactness (\cref{lemma:expdescent}).

As in the convex case, the optimality conditions of the subproblems play a vital role in the analysis. 
Since the subproblems are only solved approximately, we must again consider the inexactness of the solutions in these two steps. 
For the first step, we use strong convexity of $f$ and convexity of $g$ with appropriate perturbations to account for the inexactness. 
We call these conditions \emph{perturbed convexity} conditions, as outlined in \cref{lemma:xzconvexity}. 
\begin{lemma}[Perturbed convexity]
\label{lemma:xzconvexity}
Instate the assumptions of \cref{thm:linearconvabs}.
Let $\tilde x^{k+1}$ and $\tilde z^{k+1}$ be the inexact solutions of $x$ and $z$-subproblems under \cref{def:TypeInexactness}. Recall $(x^\star, z^\star, u^\star)$ is a saddle point of \eqref{eq:ADMMProb}. Then for some constant $C\geq 0$, the following inequalities are satisfied: 
\begin{enumerate}
\item \textup{(Semi-inexact $f$-strong convexity)}
\begin{align}
 \label{eq:lcxstrongconvexity}
   & \langle \tilde x^k-\tilde x^{k+1},\tilde x^{k+1}-x^\star\rangle_{\Theta^k} + \rho\langle N(\tilde z^{k+1}-\tilde z^{k})+u^\star-\tilde u^{k+1},M(\tilde x^{k+1}-x^{\star})\rangle + C\varepsilon_x^k \\
   & \geq \sigma_f \|\tilde x^{k+1}-x^{\star}\|^2-\langle \nabla f_2(\tilde x^{k+1})-\nabla f_2(\tilde x^k),\tilde x^{k+1}-x^\star\rangle, \nonumber
\end{align}

\item
\textup{(Semi-inexact $g$-convexity)}
\begin{equation}
\label{eq:lczconvexity1}
      \langle \tilde z^{k+1} - z^\star, \rho N^T(u^\star-\tilde u^{k+1})+\Psi^k(\tilde z^k-\tilde z^{k+1})\rangle  \ge -C\sqrt{\varepsilon_z^k}.
\end{equation}
\end{enumerate}
\end{lemma}

\ifpreprint
A detailed proof of \cref{lemma:xzconvexity} is presented in \cref{section:linearconvergenceproof}.
\fi
In \cref{lemma:xzconvexity}, we call \eqref{eq:lcxstrongconvexity} and \eqref{eq:lczconvexity1} \emph{semi-inexact} (strong) convexity because $x^\star$ ($z^\star$) is part of the exact saddle point but $\tilde x^{k+1}$ ($\tilde z^{k+1}$) is an inexact subproblem solution. 
 With \cref{lemma:xzconvexity}, we establish a descent-type inequality, which takes into inexactness.  
\begin{lemma}[Inexact sufficient descent]\label{lemma:suffdescent}
Define $\epsTot = \varepsilon_x^k+\sqrt{\varepsilon_z^k}$, then the following descent condition holds.
    \begin{align}\label{eq:suffdescent}
	&\|\tilde w^k - w^\star\|_{G^k}^2 - \|\tilde w^{k+1} - w^\star\|^2_{G^k} + C\epsTot \ge \frac{1}{2}\|\tilde w^k - \tilde w^{k+1}\|_{G^k}^2 + \frac{2\sigma_f}{\rho} \|\tilde x^{k+1} - x^\star\|^2 \\
    &-\frac{2}{\rho}\langle \nabla f_2(\tilde x^{k+1})-\nabla f_2(\tilde x^k), \tilde x^{k+1}-x^\star\rangle. \nonumber
    \end{align}
\end{lemma}
\begin{proof}
Adding the inequalities \eqref{eq:lcxstrongconvexity} and \eqref{eq:lczconvexity1} together, and using the relation $M(\tilde x^{k+1}-x^\star) = \tilde u^{k+1}-\tilde u^k+ N(z^{\star}-\tilde z^{k+1})$, we reach
\begin{align*}
    & \langle \tilde x^k-\tilde x^{k+1},\tilde x^{k+1}-x^\star\rangle_{\Theta^k}+ \langle \tilde z^k-\tilde z^{k+1},\tilde z^{k+1}-z^\star\rangle_{\Psi^k+\rho N^TN}+\rho \langle \tilde u^k-
    \tilde u^{k+1},\tilde u^{k+1}-u^\star\rangle+C\varepsilon_x^k\\
    &\geq \sigma_f \|\tilde x^{k+1}-x^{\star}\|^2-\langle \nabla f_2(\tilde x^{k+1})-\nabla f_2(\tilde x^k),\tilde x^{k+1}-x^\star\rangle-C\sqrt{\varepsilon_z^k}+\rho \langle \tilde u^k-\tilde u^{k+1},N(\tilde z^{k+1}-\tilde z^{k})\rangle \\
\end{align*}
Recalling the definitions of $\tilde w$, $\epsTot$, and $G^k$, and
using the identity $\langle a-b,\Upsilon(c-d)\rangle = 1/2\left(\|a-d\|^2_\Upsilon-\|a-c\|^2_\Upsilon\right)+1/2\left(\|c-b\|^2_\Upsilon-\|d-b\|^2_\Upsilon\right)$, we arrive at
\begin{align*}\label{eq:lcdecrease}
    &\|\tilde w^k-w^\star\|_{G^k}^2-\|\tilde w^{k+1}-w^\star\|_{G^k}^2+\epsTot\geq \frac{2\sigma_f}{\rho}\|\tilde x^{k+1}-x^\star\|^2+\|\tilde w^{k+1}-\tilde w^k\|^2_{G^k}\\ 
    &+\langle \tilde u^k-\tilde u^{k+1},N(\tilde z^{k+1}-\tilde z^{k})\rangle-\frac{2}{\rho}\langle \nabla f_2(\tilde x^{k+1})-\nabla f_2(\tilde x^k),\tilde x^{k+1}-x^\star\rangle.
\end{align*} 

Now, for the term $\|\tilde w^{k+1}-\tilde w^k\|^2_{G^k}+\langle \tilde u^k-\tilde u^{k+1},N(\tilde z^{k+1}-\tilde z^{k})\rangle$, Cauchy-Schwarz implies
\begin{align*}
    &\|\tilde w^{k+1}-\tilde w^k\|_{G^k}^2+\langle \tilde u^k-\tilde u^{k+1},N(\tilde z^{k+1}-\tilde z^k)\rangle \\
    &\geq \|\tilde w^{k+1}-\tilde w^k\|_{G^k}^2-\frac{1}{2}\|\tilde u^{k+1} -\tilde u^k\|^2-\frac{1}{2}\|\tilde z^{k+1}-\tilde z^{k}\|^2_{N^TN}\\
    &\geq \|\tilde w^{k+1}-\tilde w^k\|_{G^k}^2-\frac{1}{2}\|\tilde u^{k+1} -\tilde u^k\|^2-\frac{1}{2}\|\tilde z^{k+1}-\tilde z^{k}\|^2_{1/\rho \Psi^k+N^TN}\\
    & = \frac{1}{2}\|\tilde w^{k+1}-\tilde w^k\|_{G^k}^2.
\end{align*}
Hence we obtain
\begin{align*}
    &\|\tilde w^k-w^\star\|_{G^k}^2-\|\tilde w^{k+1}-w^\star\|_{G^k}^2+C\epsTot\geq \frac{2\sigma_f}{\rho}\|\tilde x^{k+1}-x^\star\|^2+ \frac{1}{2}\|\tilde w^{k+1}-\tilde w^k\|^2_{G^k}\\
    &-\frac{2}{\rho}\langle \nabla f_2(\tilde x^{k+1})-\nabla f_2(\tilde x^k),\tilde x^{k+1}-x^\star\rangle,
\end{align*} 
as desired. $\Box$
\end{proof}

Given the inexact sufficient descent condition \eqref{eq:suffdescent}, the next step in proving linear convergence is to 
show \eqref{eq:suffdescent} leads to a contraction relation between $\Phi^{k+1}$ and $\Phi^k$.
\begin{lemma}[Inexact Contraction Lemma]\label{lemma:expdescent}
Under the assumptions of \cref{thm:linearconvabs}, there exists constants $\delta > 0$, and $C\geq 0$ such that 
\begin{equation}\label{eq:expdescent}
(1+\delta)\Phi^{k+1}\leq (1+\zeta^k)\left(\Phi^k+C\epsTot\right).
\end{equation}
\end{lemma}
\ifpreprint
The proof of \cref{lemma:expdescent}, and an explicit expression for $\delta$, appears in \cref{section:linearconvergenceproof}.
\fi
As in \cref{thm:linearconvabs}, the constant $\delta$ gives the rate of linear convergence and depends on the conditioning of $f$ and the constraint matrices $M$ and $N$.
The better the conditioning, the faster the convergence, with the opposite holding true as the conditioning worsens. 
With \cref{lemma:expdescent} in hand, we now prove \cref{thm:linearconvabs}. 

\subsection{Proof of \cref{thm:linearconvabs}}
\label{subsection:LinConPf}
\proof
By induction on \eqref{eq:expdescent}, we have 
\begin{align*}\label{eq:inexactsumabs}
(1+\delta)^k\Phi^k & \leq \left(\prod^{k}_{j=1}(1+\zeta^j)\right)\Phi^1+C\sum^{k}_{j=1}\left(\prod^{k}_{i=j}(1+\zeta^{i})\right)(1+\delta)^{j-1}\epsjTot \\
&\leq \tau_\zeta \Phi^1+C\tau_\zeta\sum_{j=1}^{k}(1+\delta)^{j-1}\epsjTot \leq \tau_\zeta \Phi^1+C\tau_\zeta \sum_{j=1}^{k}\left(\frac{1+\delta}{1+q}\right)^{j-1} \\
&\leq \tau_\zeta\Phi^1+ \frac{C\tau_\zeta}{1-\frac{1+\delta}{1+q}} = \tau_\zeta\Phi^1+S,
\end{align*}
where the second inequality uses Assumption \ref{assp:geodecay} to reach $\epsjTot \leq \frac{C}{(1+q)^{j}}$, and the third inequality uses $q>\delta$ to bound the sum by the sum of the geometric series.
Hence, we have shown the first claim.
The second claim follows immediately from the first via a routine calculation. $\Box$

\section{Applications}
\label{section:applications}
This section applies our theory to establish convergence rates for NysADMM and sketch-and-solve ADMM.

\subsection{Convergence of NysADMM}
\label{subsection:NysADMM}
\begin{algorithm}[H]
	\centering
	\caption{NysADMM \label{alg-nysadmm}}
	\begin{algorithmic}
		\Require{penalty parameter $\rho$, step-size $\eta$, regularization $\sigma \geq 0$, forcing sequences $\{\varepsilon_x^k\}_{k\geq 1}$, $\{\varepsilon_z^k\}_{k\geq 1}$}
		\Repeat
        \State{Find $\varepsilon_x^k$-approximate solution $\tilde{x}^{k+1}$ of 
        \[(\eta H_f^k+(\rho+\eta\sigma) I)x =  \eta (H^k_f+\sigma I)\tilde x^k - \nabla f(\tilde x^k)+\rho(\tilde z^k - \tilde u^k)\]}
		\State{$\tilde{z}^{k+1} \overset{\varepsilon_z^k}{\approxeq} \underset{z\in Z}{\argmin} \{g(z) + \frac{\rho}{2}\|M\tilde x^{k+1} -z + \tilde u^k\|^2\}$}
		\State{$\tilde u^{k+1} = \tilde u^k + M\tilde x^{k+1} -\tilde z^{k+1}$}
		\Until{convergence}
		\Ensure{solution $(x^\star, z^\star)$ of problem \eqref{eq:ADMMProb}}
	\end{algorithmic}
\end{algorithm}
We begin with the NysADMM scheme from \cite{pmlr-v162-zhao22a}. 
Recall NysADMM is obtained from \cref{alg-gen_newton_admm} by setting $f_1 = 0, f_2 = f$,
using the exact Hessian $\Theta^k = H_f(\tilde x^k) = \eta (H_f^k+\sigma I)$, and setting $\Psi^k = 0$.
Instantiating these selections into \cref{alg-gen_newton_admm},
we obtain NysADMM, presented as \cref{alg-nysadmm}.
Compared to the original NysADMM, \cref{alg-nysadmm} adds a regularization term $\sigma I$ to the Hessian. 
In theory, when $f$ is only convex, 
this regularization term is required to ensure the condition $\Theta^k \succeq \theMin$ of Assumption \ref{assp:regularity}, 
as the Hessian along the optimization path may fail to be uniformly bounded below.
The addition of the $\sigma I$ term removes this issue. 
However, the need for this term seems to be an artifact of the proof, 
as \citet{pmlr-v162-zhao22a} runs \cref{alg-nysadmm} on non-strongly convex objectives with $\sigma = 0$, and convergence is still obtained.
Similarly, 
we set $\sigma = 0$ in all our experiments, and convergence consistent with our theory is obtained. 
Hence, in practice, 
we recommend setting $\sigma = 0$, or some small value, say $\sigma = 10^{-8}$, 
so that convergence isn't slowed by unneeded regularization. 

The $x$-subproblem for NysADMM is a linear system. 
NysADMM solves this system using the Nystr{\"o}m PCG method from \cite{frangella2023randomized}.
The general convergence of NysADMM was left open in \cite{pmlr-v162-zhao22a},
which established convergence only for quadratic $f$.
Moreover, the result in \cite{pmlr-v162-zhao22a} does not provide an explicit convergence rate.
We shall now rectify this state of affairs using \cref{thm:generalconv}.
We obtain the following convergence guarantee by substituting the parameters defining NysADMM (with the added $\sigma I$ term) into \cref{thm:generalconv}.
\begin{corollary}[Convergence of NysADMM]
\label{thm:NysADMM}
Instate the assumptions of \cref{thm:generalconv}. 
Let $\sigma >0$. Set $f_1 = 0, f_2 = f$, $\eta = \hat{L}_f$, $\Theta^k = \eta (H^k_f+\sigma I)$, and $\Psi^k = 0$ in \cref{alg-gen_newton_admm}. 
Then
\begin{align*}
	& f(\bar x^{t+1})+g(\bar z^{t+1})-p^{\star} \leq \frac{\Gamma}{t},
        & \|M\bar x^{t+1}+N\bar z^{t+1}\| \le \frac{2}{t}\sqrt{\frac{\Gamma}{\rho}+d^2_{u^\star}}.
\end{align*}
Here $\Gamma$ and $d_{u^\star}^2$ are the same as in \cref{thm:generalconv}. 
\end{corollary}
NysADMM converges at the same $\mathcal O(1/t)$-rate as standard ADMM, 
despite all the approximations it makes. 
Thus, NysADMM offers the same level of performance as ADMM, 
but is much faster due to its use of inexactness.
This result is empirically verified in \cref{section:experiments}, where NysADMM converges almost identically to ADMM.
\cref{thm:NysADMM} supports the empirical choice 
of a constant step-size $\eta = 1$,
which was shown to have excellent performance uniformly across tasks in \cite{pmlr-v162-zhao22a}:
the theorem sets $\eta = \hat L_f$
and $\hat L_f = 1$ for quadratic functions,
and satisfies $\hat L_f = \mathcal O(1)$ for loss functions such as the logistic loss.
We recommend setting $\eta = 1$ as the default value for GeNI-ADMM. 
Given NysADMM's superb empirical performance in \cite{pmlr-v162-zhao22a} 
and the firm theoretical grounding given by \Cref{thm:NysADMM}, 
we conclude that NysADMM provides a reliable framework 
for solving large-scale machine learning problems.

\subsection{Convergence of sketch-and-solve ADMM}
\label{subsection:SketchAndSolveADMM}
\begin{algorithm}[H]
	\centering
	\caption{Sketch-and-solve ADMM \label{alg-sketch-solve-admm}}
	\begin{algorithmic}
		\Require{penalty parameter $\rho$, step-size $\eta$, $\{\varepsilon_z^k\}_{k\geq 1}$}
		\Repeat
        \State{Find solution $\tilde{x}^{k+1}$ of 
        \[(\eta \hat H^k+(\rho+\eta\gamma^k)I)x = \eta \left(\hat H^k+\gamma^k I\right)\tilde x^k - \nabla f(\tilde x^k)+\rho(\tilde z^k - \tilde u^k)\]}
		\State{$\tilde{z}^{k+1} \overset{\varepsilon_z^k}{\approxeq} \underset{z\in Z}{\argmin} \{g(z) + \frac{\rho}{2}\|M\tilde x^{k+1} -z + \tilde u^k\|^2\}$}
		\State{$\tilde u^{k+1} = \tilde u^k + M\tilde x^{k+1} -\tilde z^{k+1}$}
		\Until{convergence}
		\Ensure{solution $(x^\star, z^\star)$ of problem \eqref{eq:ADMMProb}}
	\end{algorithmic}
\end{algorithm}
Sketch-and-solve ADMM is obtained from GeNI-ADMM by setting $\Theta^k$ 
to be an approximate Hessian computed by a sketching procedure.
The two most popular sketch-and-solve methods are the Newton sketch \citep{pilanci2017newton,lacotte2021adaptive,derezinski2021newton} 
and Nystr{\"o}m sketch-and-solve \citep{bach2013sharp,alaoui2015fast,frangella2023randomized}.
Both methods use an approximate Hessian that is cheap to compute
and yields a linear system that is fast to solve.
ADMM together with sketching techniques, provides a compelling means
to handle massive problem instances.
The recent survey \cite{buluc2021randomized} suggests using sketching to solve a quadratic ADMM subproblem for efficiently solving large-scale inverse problems.
However, it was previously unknown whether the resulting method would converge.
Here, we formally describe the sketch-and-solve ADMM method and prove convergence. 

Sketch-and-solve ADMM is obtained from \cref{alg-gen_newton_admm} by setting 
\begin{equation}
\label{eq:sketch-and-solve}
\Theta^k = \eta\left(\hat{H}^k+\gamma^k I\right),
\end{equation}
where $\hat{H}^k$ is an approximation to the Hessian $H_f(\tilde x^k)$ at the $k$th iteration, 
and $\gamma^k \geq 0$ is a constant chosen to ensure convergence. 
The term $\gamma^k I$ ensures that the approximate linearization 
satisfies the $\Theta$-relative smoothness condition
when $\gamma^k$ is chosen appropriately, as in the following lemma:
\begin{lemma}
\label{lemma:RelSmoothSS}
    Suppose $f$ is $\hat{L}_f$-relatively smooth with respect to its Hessian $H_f$. Construct $\{\Theta^k\}_{k\ge 1}$ as in \eqref{eq:sketch-and-solve} and select $\gamma^k>0$ such that $\gamma^k\geq \|E^k\| = \|H_f(\tilde x^k)-\hat{H}^k\|$ for every $k$. 
    Then
    \[f(x)\leq f(\tilde x^k)+\langle \nabla f(\tilde x^k),x-\tilde x^k\rangle +\frac{\hat{L}_f}{2}\|x-\tilde x^k\|^2_{\Theta^k}.\]
\end{lemma}

\cref{lemma:RelSmoothSS} shows we can ensure relative smoothness by selecting $\gamma^k>0$ appropriately. 
Assuming relative smoothness, 
we may invoke \cref{thm:generalconv} to guarantee that sketch-and-solve ADMM converges.
Unlike with NysADMM, we find it is necessary to select $\gamma^k$ carefully (such as in \cref{lemma:RelSmoothSS}) to ensure the relative smoothness condition holds, 
otherwise sketch-and-solve ADMM will diverge, 
see \cref{subsection:correction_exp} for numerical demonstration.
This condition is somewhat different from those in prior sketching schemes in optimization, such as the Newton Sketch \citep{pilanci2017newton,lacotte2021adaptive}, where convergence is guaranteed as long as the Hessian approximation is invertible. 

\begin{corollary}
\label{thm:sketch-solve}
Suppose $f$ is $\hat{L}_f$-relatively smooth with respect to its Hessian $H_f$ and instate the assumptions of \cref{thm:generalconv}. 
In \cref{alg-gen_newton_admm}, 
Set $\Theta^k =  \eta(\hat H^k +\gamma^k I)$ with $\eta = \hat L_f$, $\gamma^k = \|E^k\| = \|\hat H^k-H_f(\tilde x^k)\|$, and $\Psi^k = 0$.
Then
\begin{equation*}
	f(\bar x^{t+1})+g(\bar z^{t+1})-p^{\star} \leq \frac{1}{t}\left(\hat{L}_f \left(d^2_{x^\star,H_f(\tilde x^1)}+\|E^1\|d^2_{x^\star}\right) + \frac{\rho}{2}d^2_{z^\star} + C_x\mathcal{E}_x + C_z\mathcal{E}_z + C_\zeta \mathcal{E}_\zeta\right)=:\frac{1}{t}\hat \Gamma
	\end{equation*}
    and
    \begin{equation*}
        \|M\bar x^{t+1} - \bar z^{t+1}\| \le \frac{2}{t}\sqrt{\frac{\hat \Gamma}{\rho} + d^2_{u^\star}}.
    \end{equation*}
Here variables $\bar x^{t+1}$ and $\bar z^{t+1}$, diameters $d^2_{x^\star,H_f(\tilde x^1)}$, $d^2_{z^\star}$, and $d^2_{u^\star}$, and constants $C_x$, $C_z$, $\mathcal E_x$, $\mathcal E_z$, $\mathcal E_\zeta$, $R_\star$, and $p^\star$ are all defined as in \cref{thm:generalconv}.
\end{corollary}

\begin{remark}
	It is not hard to see \cref{thm:sketch-solve} also holds 
	for any $\gamma>0$ if the $\{\Theta^{k}\}_{k\ge 1}$ satisfy 
	\[
	(1-\tau)(H_f^k+\gamma I)\preceq \Theta^k\preceq (1+\tau) (H_f^k+\gamma I),
	\]
	for some $\tau\in (0,1)$, 
	and $\eta \geq \hat{L}_f/(1-\tau)$. 
	Many sketching methods can ensure this relative error condition with high probability \citep{Karimireddy2018global}
	by choosing the sketch size proportional to the effective dimension of 
	$H_f + \gamma I$ \citep{frangella2023randomized,lacotte2021adaptive}.
\end{remark}  

\cref{thm:sketch-solve} shows sketch-and-solve ADMM obtains an $\mathcal O(1/t)$-convergence rate. 
The main difference between NysADMM and sketch-and-solve is 
the additional error term due to the use of the approximate Hessian,
which results in a slightly slower convergence rate. 
In this sense, sketch-and-solve ADMM can be regarded as a compromise between NysADMM ($\Theta^k = \eta (H_f(\tilde x^k)+\sigma I)$) and gradient descent ADMM ($\Theta^k = \eta I$) --- with the convergence rate improving as the accuracy of the Hessian approximation increases.  
\section{Numerical experiments}
\label{section:experiments}
In this section, we numerically illustrate the convergence results developed in \cref{section:SubLinConvergence} 
for several methods highlighted in \cref{subsection:recoveredadmm} that fit into the GeNI-ADMM framework: 
sketch-and-solve ADMM (\cref{alg-sketch-solve-admm}),
NysADMM (\cref{alg-nysadmm}), and ``gradient descent'' ADMM (GD-ADMM) \eqref{eq:gd-admm-update}.
As a baseline, we also compare to exact ADMM (\cref{alg-admm}) to see how various approximations or inexactness impact the convergence rate.
We conduct three sets of experiments that verify different aspects of the theory:
\begin{itemize}
    \item \cref{subsection:sublin_exp} verifies that NysADMM, GD-ADMM and sketch-and-solve ADMM converge sublinearly for convex problems.
    Moreover, we observe a fast transition to linear convergence, 
    after which all methods but GD-ADMM converge quickly to high accuracy. 
    \ifpreprint
    \item \cref{subsection:lin_exp} verifies that NysADMM, GD-ADMM, and sketch-and-solve ADMM converge linearly for strongly convex problems. 
    \else
    \fi
    \item \cref{subsection:correction_exp} verifies that, without the correction term, sketch-and-solve ADMM diverges, showing the necessity  of the correction term in \cref{subsection:SketchAndSolveADMM}.  
\end{itemize}

We consider three common problems in machine learning and statistics in our experiments: lasso \citep{tibshirani1996regression}, elastic net regression \citep{zou2005regularization}, and $\ell_1$-logistic regression \citep{hastie2015statistical}.
All experiments use the \texttt{realsim} dataset from LIBSVM~\citep{libsvm}, accessed through OpenML~\citep{openml}, which has $72,309$ samples and $20,958$ features.
Our experiments use a subsample of \texttt{realsim}, consisting of $10,000$ random samples, which ensures the objective is not strongly convex for lasso and $\ell_1$-logistic regression.


For NysADMM, a sketch size of $50$ is used to construct the Nystr{\"o}m preconditioner, and for sketch-and-solve ADMM, we use a sketch size $500$ to form the Hessian approximation.
For sketch-and-solve ADMM, the parameter $\gamma^k$ in~\eqref{eq:sketch-and-solve} is chosen by estimating the error of the Nystr\"{o}m sketch using power iteration.


All experiments are performed in the Julia programming language~\citep{bezanson2017julia} on a 
MacBook Pro with a M1 Max processor and 64GB of RAM.
To compute the ``true'' optimal values, we use the commercial solver Mosek~\citep{mosek} (with tolerances set low for high accuracy and presolve turned off to preserve the problem scaling) and the modeling language JuMP~\citep{dunning2017jump,legat2021mathoptinterface}.

\subsection{GeNI-ADMM converges sublinearly and locally linearly on convex problems}
\label{subsection:sublin_exp}
To illustrate the global sublinear convergence of GeNI-ADMM methods on convex objectives, we look at the performance of ADMM, NysADMM, GD-ADMM, 
and sketch-and-solve ADMM on solving a lasso and $\ell_1$-logistic regression problem with the \texttt{realsim} dataset.
Note, as the number of samples is smaller than the number of features, 
the corresponding optimization problems are convex but not strongly convex. 
\paragraph{Lasso regression}
\begin{figure}[t]
\captionsetup[sub]{font=scriptsize}
    \centering
    \begin{subfigure}[t]{0.4\textwidth}
        \centering
        \includegraphics[width=\columnwidth]{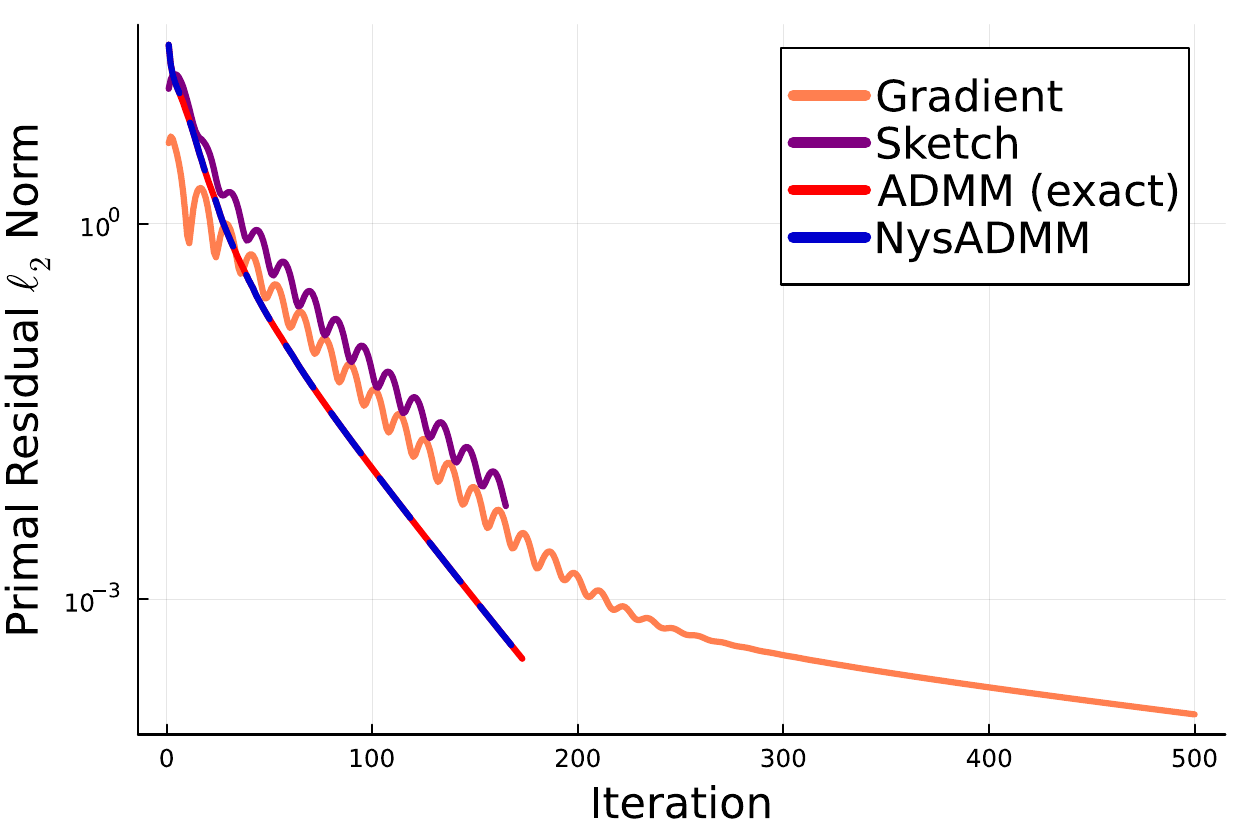}
        \caption{Primal residual convergence}
    \end{subfigure}
    \hfill
    \begin{subfigure}[t]{0.4\textwidth}
        \centering
        \includegraphics[width=\columnwidth]{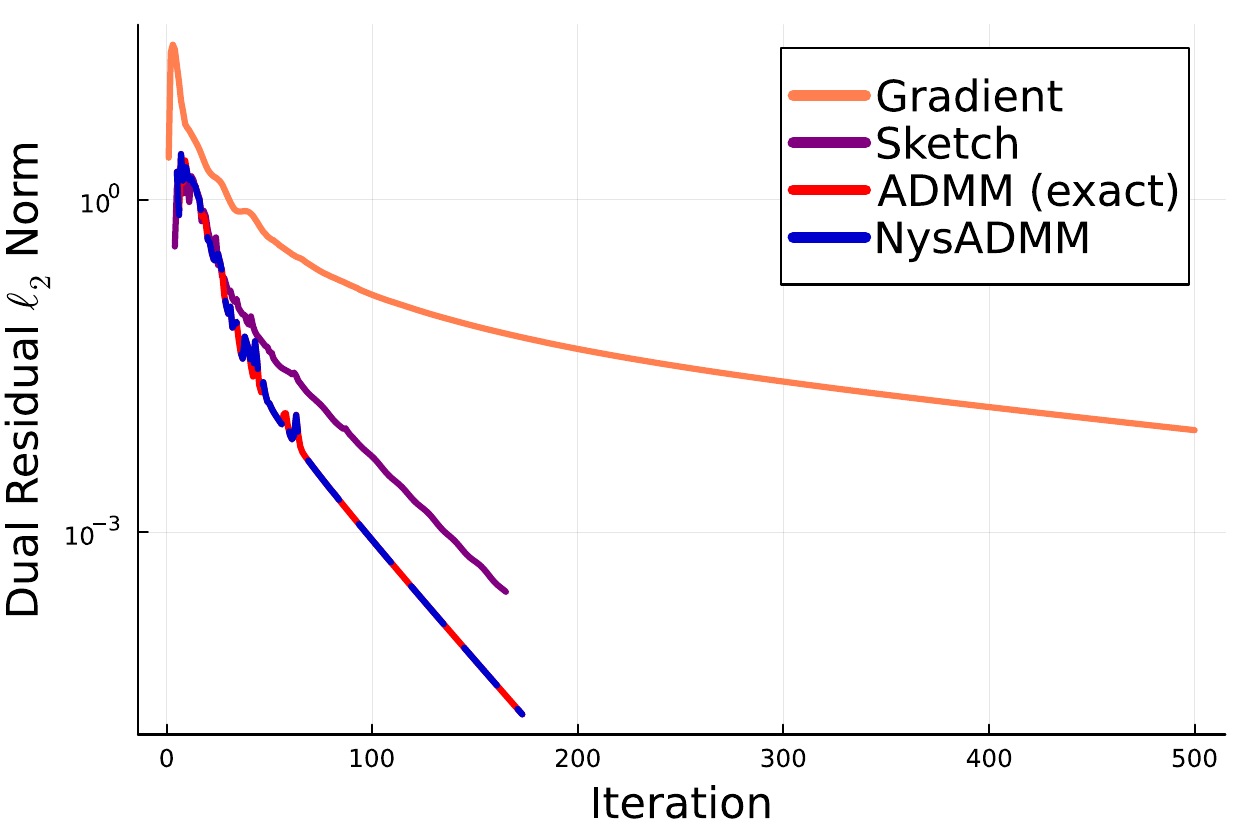}
        \caption{Dual residual convergence}
    \end{subfigure}
    \vskip\baselineskip
    \begin{subfigure}[t]{0.4\textwidth}
        \centering
        \includegraphics[width=\columnwidth]{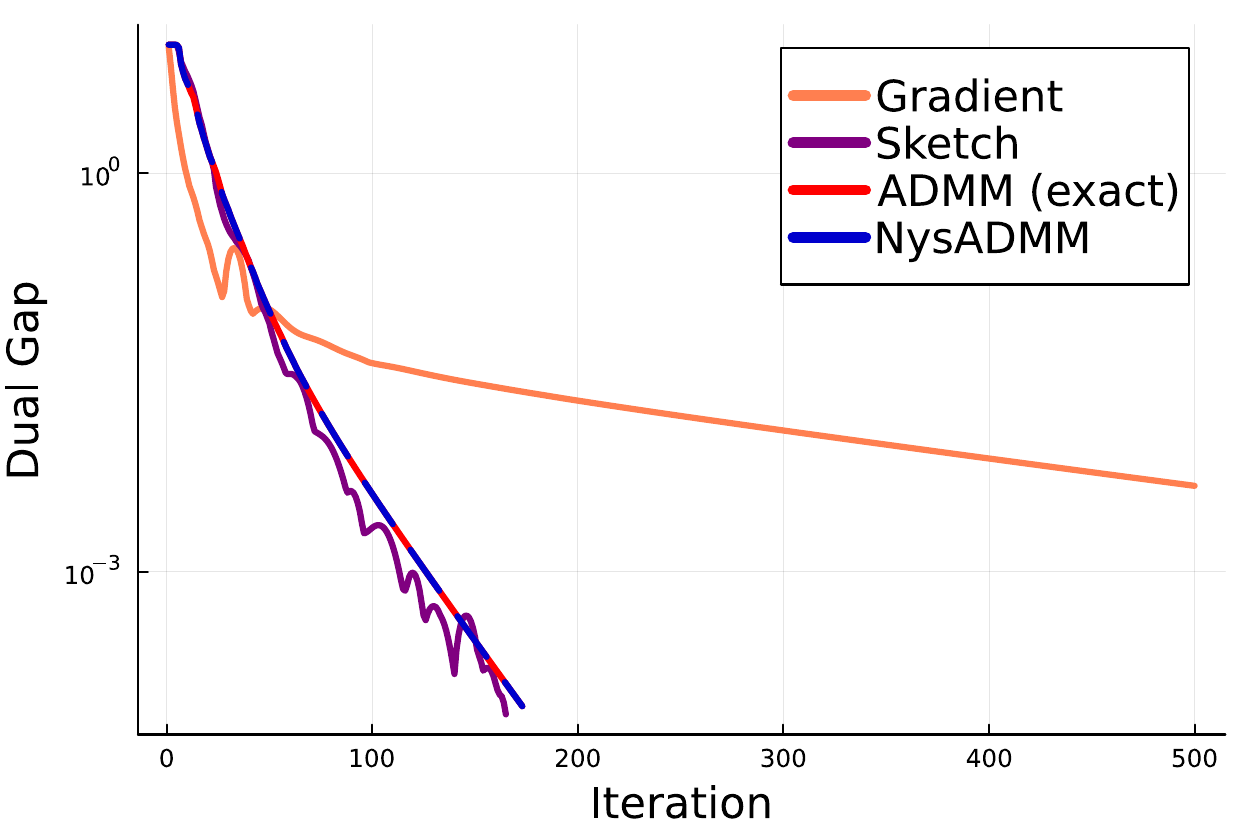}
        \caption{Dual gap convergence}
    \end{subfigure}
    \hfill
    \begin{subfigure}[t]{0.4\textwidth}
        \centering
        \includegraphics[width=\columnwidth]{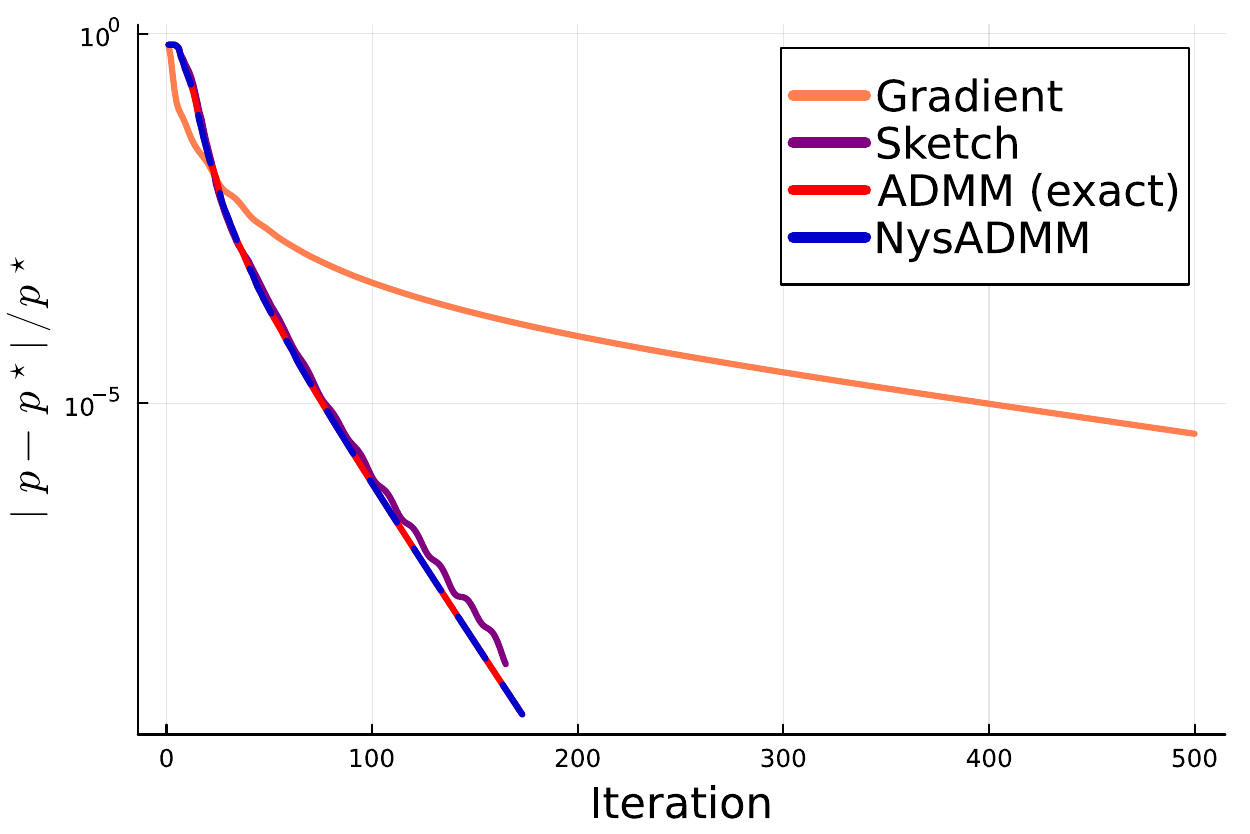}
        \caption{Objective value convergence}
    \end{subfigure}
    \caption{Convergence of lasso regression for NysADMM, sketch-and-solve ADMM, and gradient descent ADMM.}
    \label{fig:lasso}
\end{figure}
The lasso regression problem is to minimize the $\ell_2$ error of a linear model with an $\ell_1$ penalty on the weights:
\[
\begin{aligned}
& \text{minimize} && ({1}/{2})\|Ax - b\|_2^2 + \gamma \|x\|_1.
\end{aligned}
\]
This can be easily transformed into the form~\eqref{eq:ADMMProb} by taking $f(x) = (1/2)\|Ax - b\|_2^2$, $g(z) = \gamma \|z\|_1$,  $M = I$, $N = -I$,
and $X = Z = \reals^n$. 
We set $\gamma = 0.05 \cdot \gamma_\mathrm{max}$, where $\gamma_\mathrm{max} = \|A^Tb\|_\infty$ is the value above which
the all zeros vector is optimal \cite{hastie2015statistical}.
We stop the algorithm when the gap is less than $10^{-4}$, or after $500$ iterations.

The results of lasso regression are illustrated in \Cref{fig:lasso}.
\Cref{fig:lasso} shows all methods initially converge at sublinear rate, but quickly transition to a linear rate of convergence after reaching some manifold containing the solution. 
ADMM, NysADMM, and sketch-and-solve ADMM (which use curvature information) all converge converge much faster then GD-ADMM, confirming the predictions of \cref{section:SubLinConvergence} that methods which use curvature information will converge faster than methods that do not. 
Moreover, the difference in convergence between NysADMM and ADMM is negligible, despite the former having a much cheaper iteration complexity due to the use of inexact linear system solves. 
\paragraph{L1-Logistic Regression}
\begin{figure}[t] 
\captionsetup[sub]{font=scriptsize}
    \centering
    \begin{subfigure}[t]{0.4\textwidth}
        \centering
        \includegraphics[width=\columnwidth]{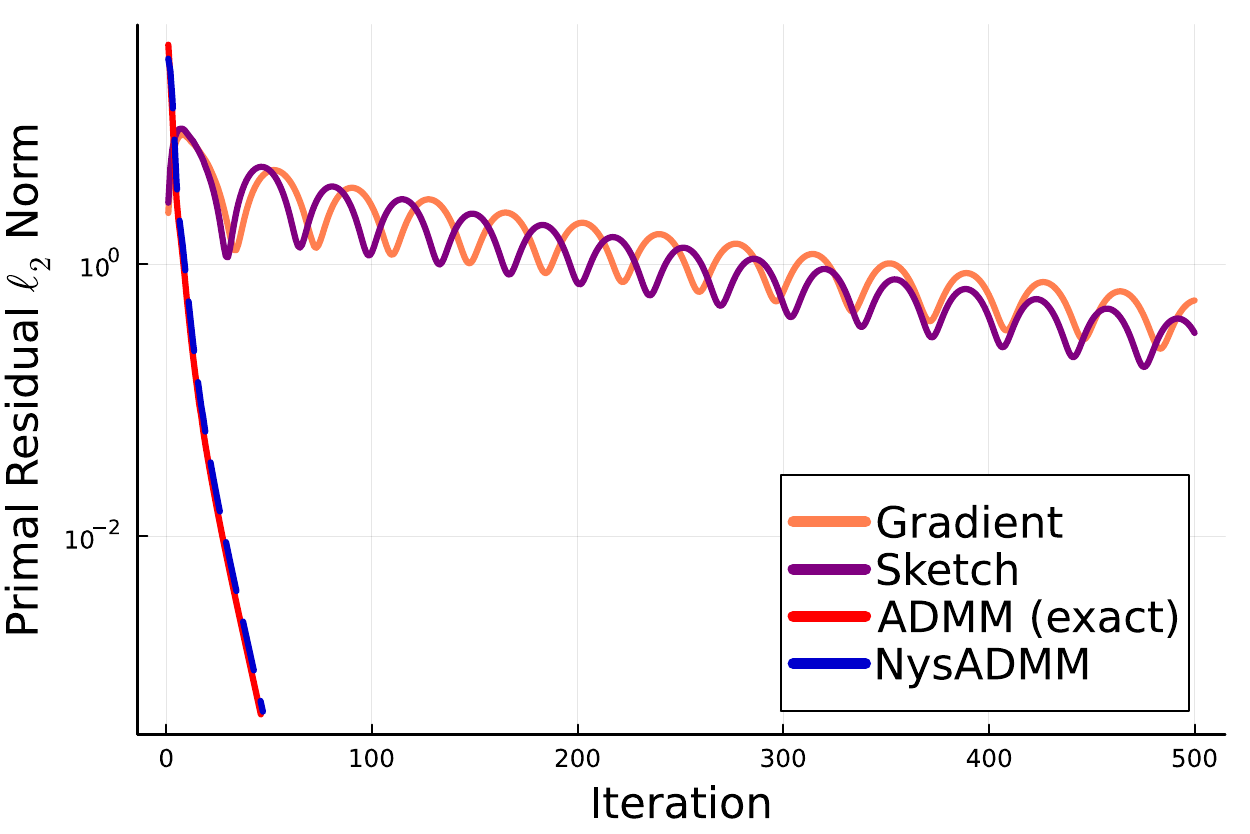}
        \caption{Primal residual convergence}
    \end{subfigure}
    \hfill
    \begin{subfigure}[t]{0.4\textwidth}
        \centering
        \includegraphics[width=\columnwidth]{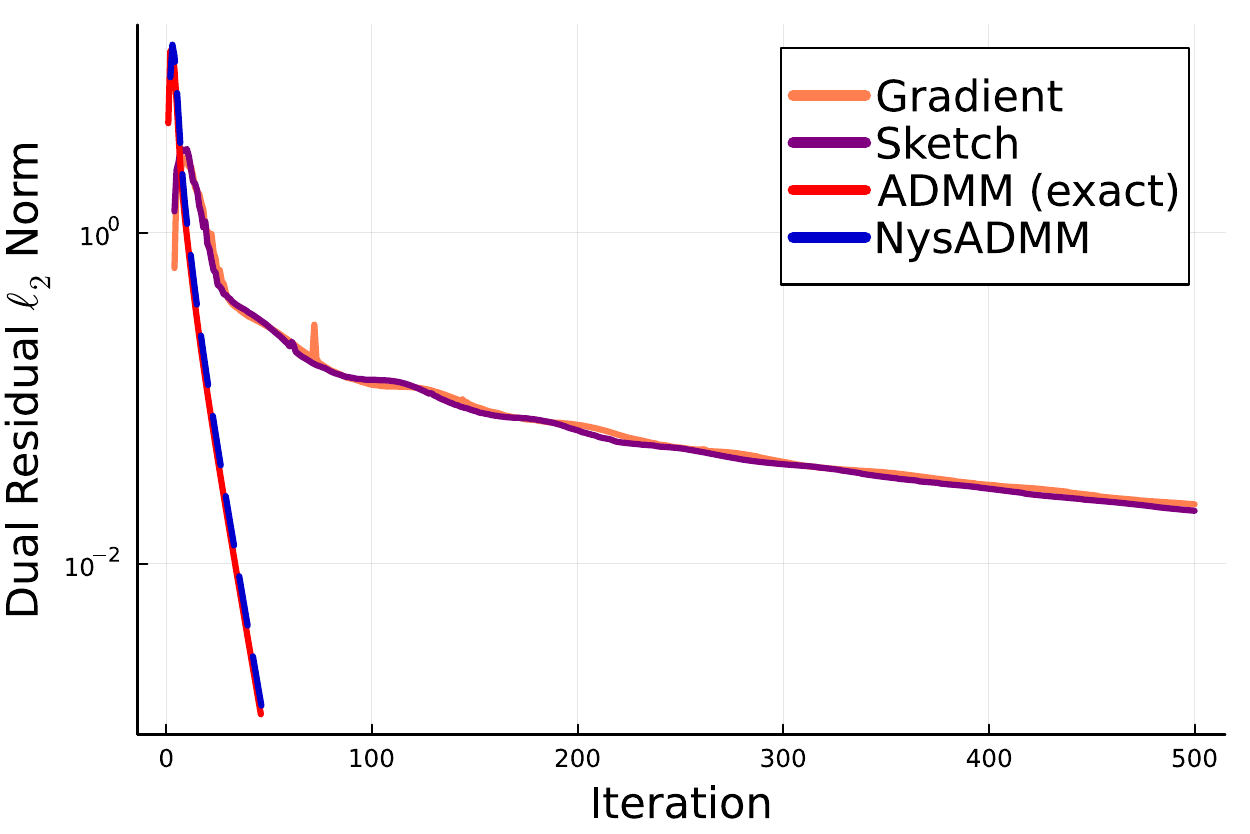}
        \caption{Dual residual convergence}
    \end{subfigure}
    \vskip\baselineskip
    \begin{subfigure}[t]{0.4\textwidth}
        \centering
        \includegraphics[width=\columnwidth]{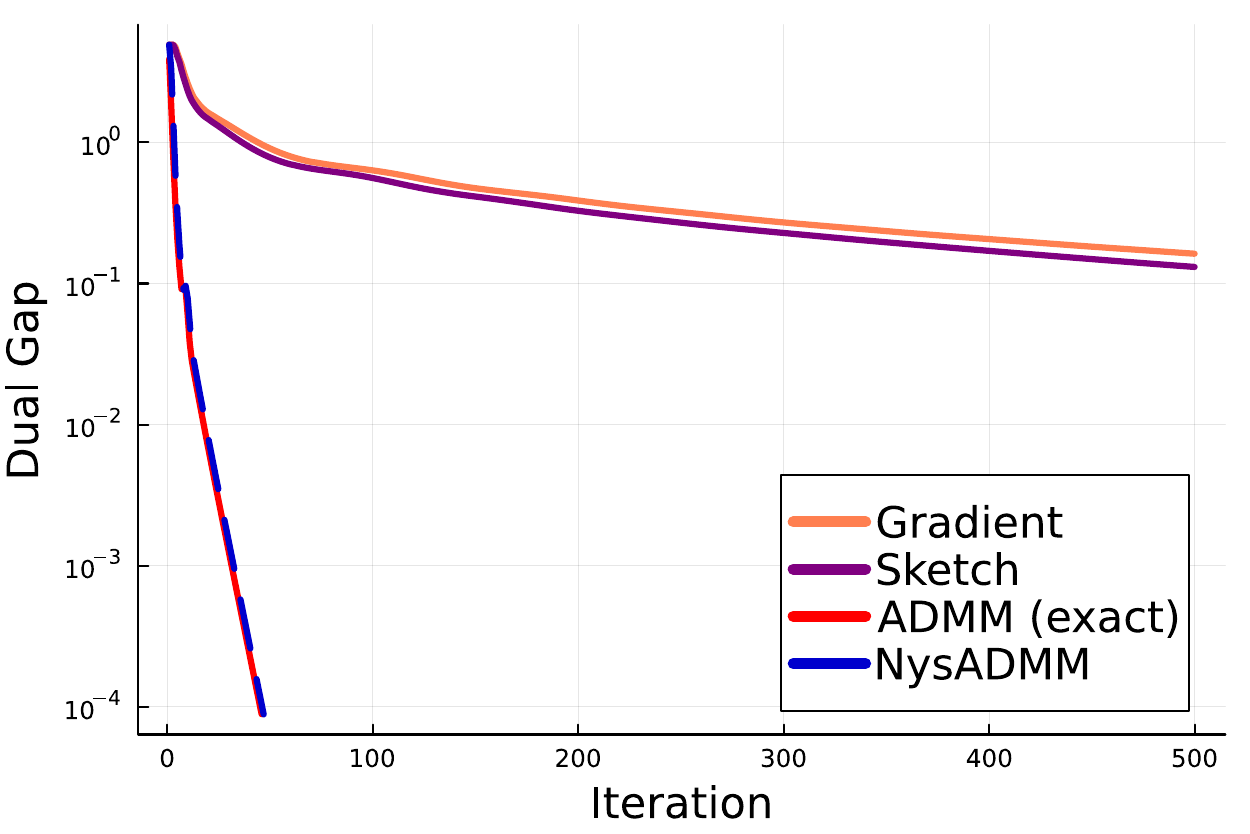}
        \caption{Dual gap convergence}
    \end{subfigure}
    \hfill
    \begin{subfigure}[t]{0.4\textwidth}
        \centering
        \includegraphics[width=\columnwidth]{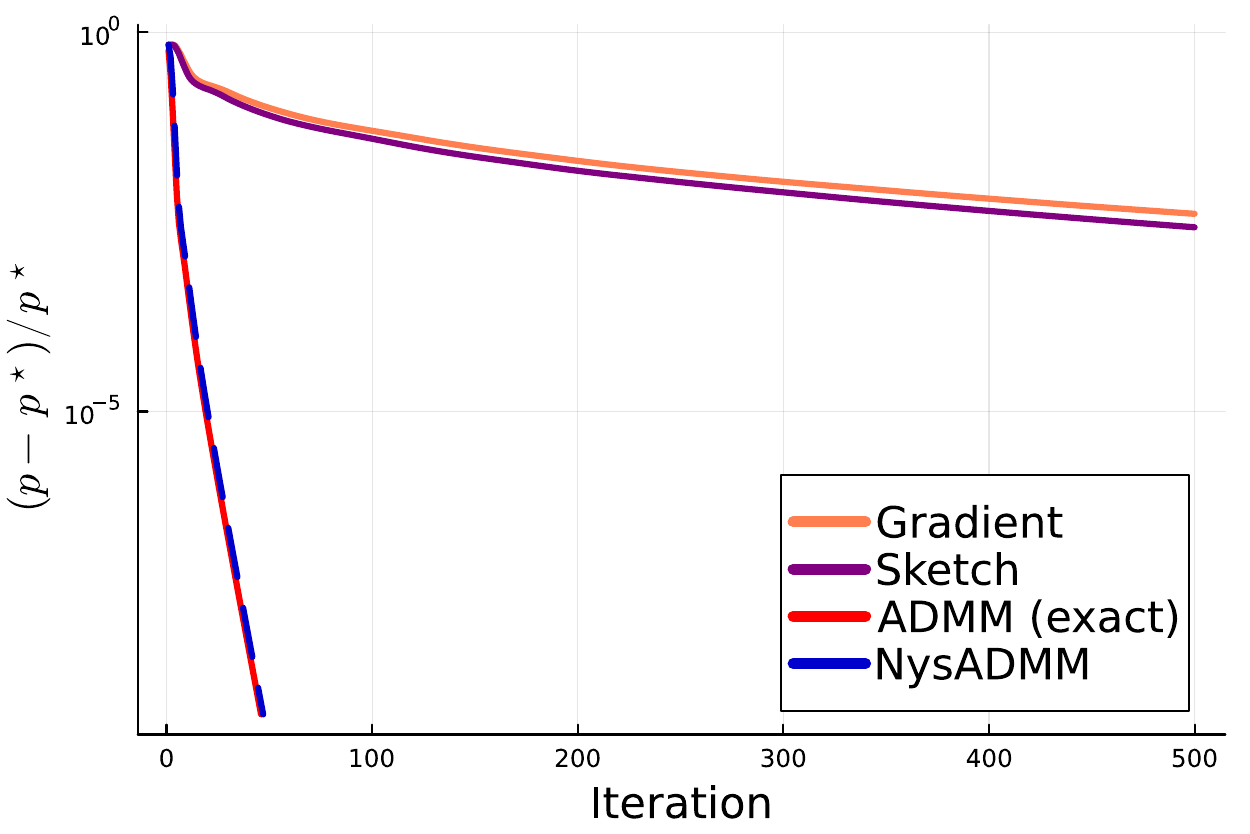}
        \caption{Objective value convergence}
    \end{subfigure}
    \caption{Convergence of logistic regression for NysADMM, sketch-and-solve ADMM, and gradient descent ADMM.}
    \label{fig:logistic}
\end{figure}

We set $\gamma = 0.05 \cdot \gamma_\mathrm{max}$, where $\gamma_\mathrm{max} = (1/2)\|A^T\mathbf{1}\|_\infty$ is the value above which
the all zeros vector is optimal. 
For NysADMM, the preconditioner is re-constructed after every 20 iterations.
For sketch-and-solve ADMM, we re-construct the approximate Hessian at every iteration.
Since the $x$-subproblem is not a quadratic program, we use the L-BFGS~\citep{liu1989limited} implementation from the \texttt{Optim.jl} package~\citep{mogensen2018optim} to solve the $x$-subproblem.

\Cref{fig:logistic} presents the results for logistic regression. 
The results are consistent with the lasso experiment---all methods initially converge sublinearly before quickly transitioning to linear convergence, and  methods using better curvature information, converge faster.
In particular, although sketch-and-solve-ADMM converges slightly faster than GD-ADMM, its convergence is much slower than NysADMM, which more accurately captures the curvature due to using the exact Hessian.
The convergence of NysADMM and ADMM is essentially identical, despite the former having a much cheaper iteration cost due to approximating the $x$-subproblem and using inexact linear system solves.

\ifpreprint
\subsection{GeNI-ADMM converges linearly on strongly convex problems}
\label{subsection:lin_exp}
\begin{figure}[t]
\captionsetup[sub]{font=scriptsize}
    \centering
    \begin{subfigure}[t]{0.46\textwidth}
        \centering
        \includegraphics[width=\columnwidth]{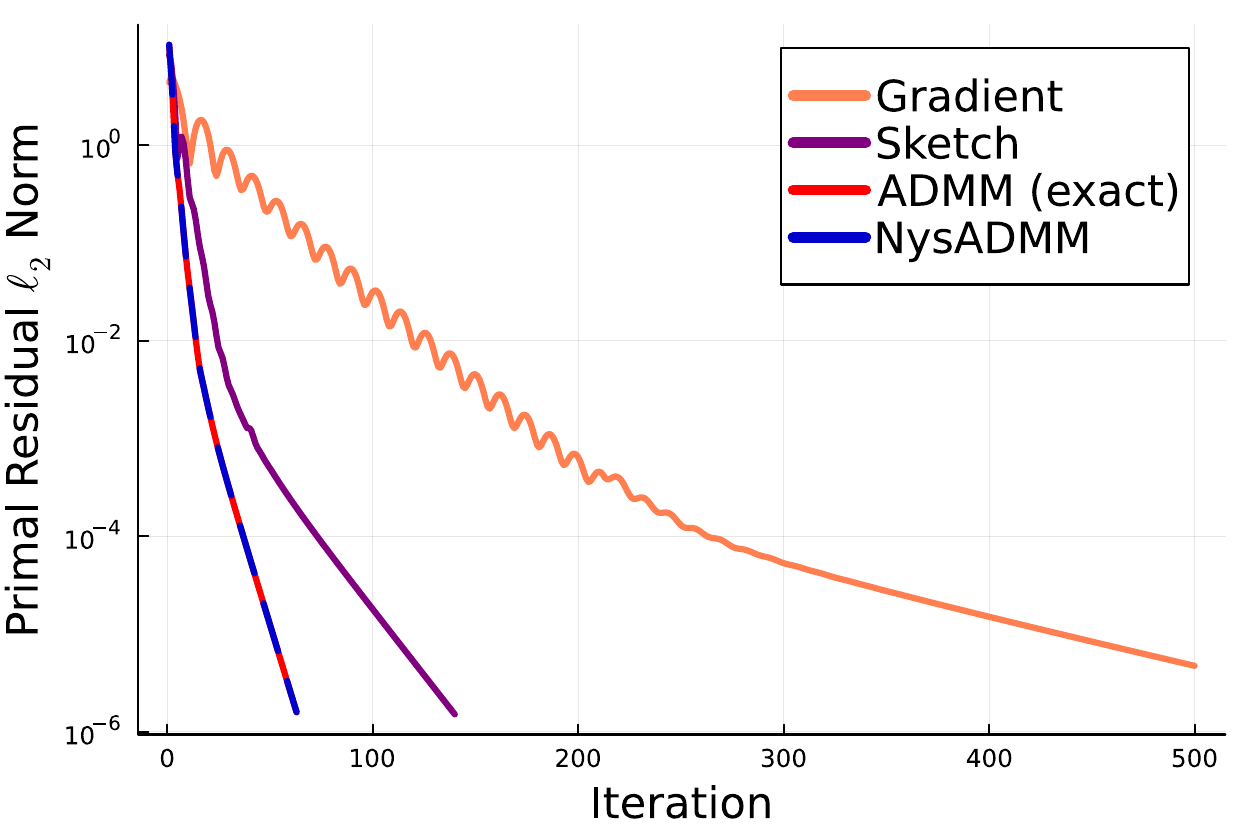}
        \caption{Primal residual convergence}
    \end{subfigure}
    \hfill
    \begin{subfigure}[t]{0.46\textwidth}
        \centering
        \includegraphics[width=\columnwidth]{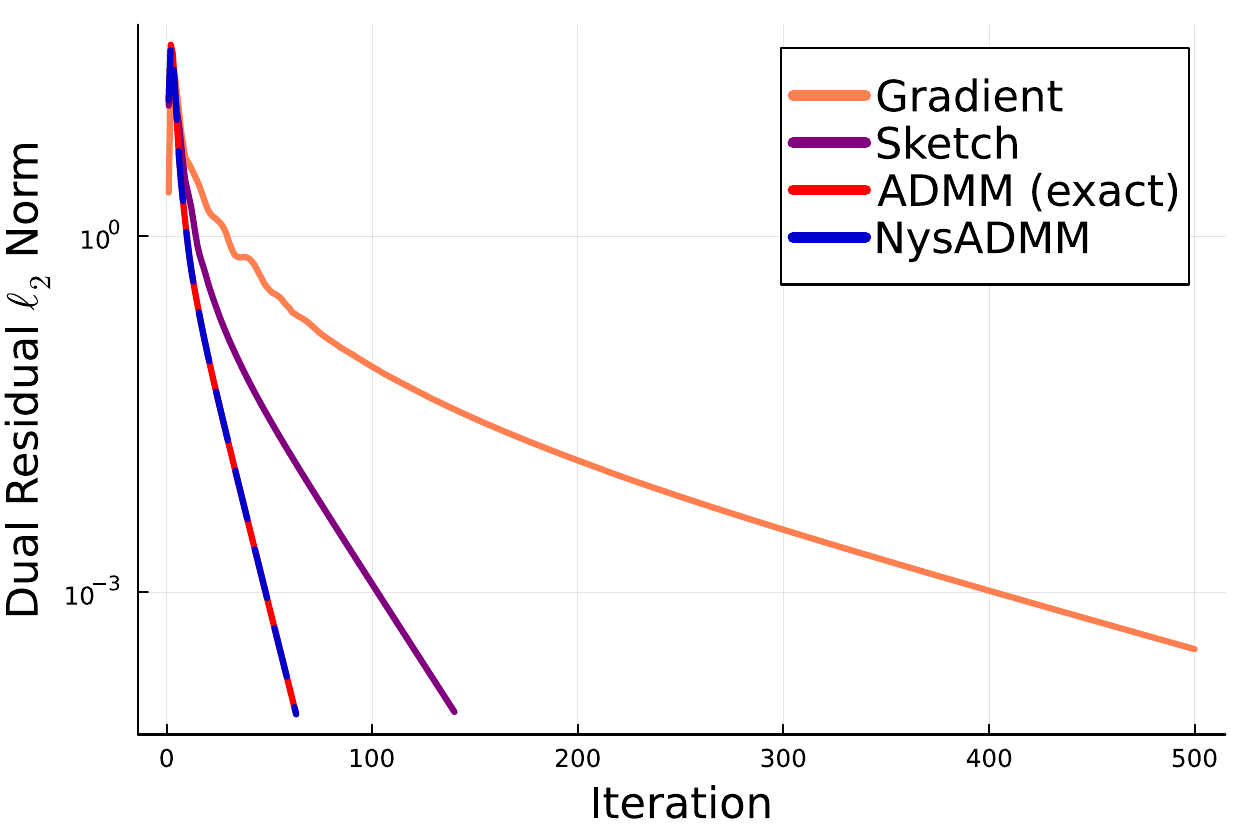}
        \caption{Dual residual convergence}
    \end{subfigure}
    \vskip\baselineskip
    \begin{subfigure}{0.46\textwidth}
        \centering
        \includegraphics[width=\columnwidth]{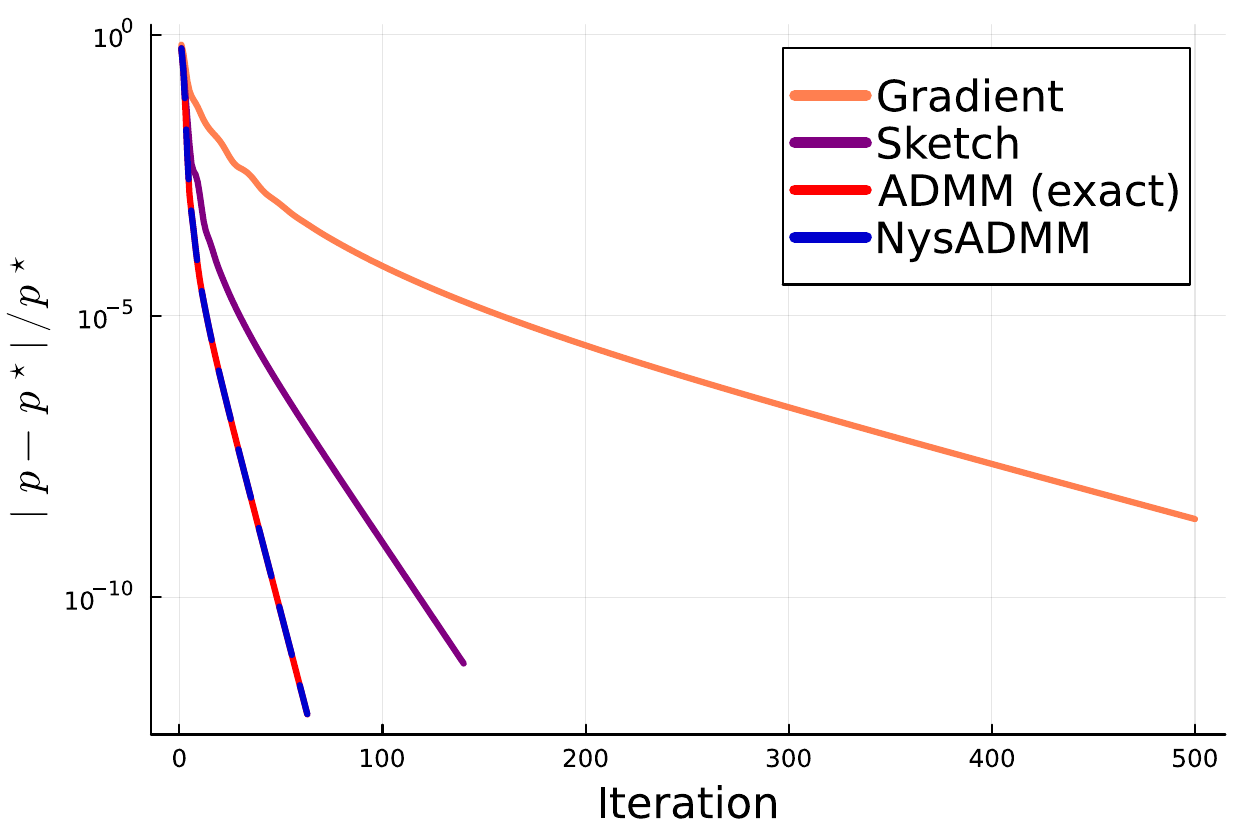}
        \caption{Objective value convergence}
    \end{subfigure}
    \caption{Convergence of elastic net regression for NysADMM, sketch-and-solve ADMM, and gradient descent ADMM.}
    \label{fig:elastic-net}
\end{figure}

To verify the linear convergence of GeNI-ADMM methods in the presence of strong convexity, we experiment with the elastic net problem:
\[
\begin{aligned}
& \text{minimize} && ({1}/{2})\|Ax - b\|_2^2 + \gamma \|x\|_1 + (\mu/2)\|x\|_2^2.
\end{aligned}
\]
We set $\mu = 1$, and use the same problem data and value of $\gamma$ as the lasso experiment.
The results of the elastic-net experiment are presented in \Cref{fig:elastic-net}.
Comparing Figures \ref{fig:lasso} and \ref{fig:elastic-net}, 
we clearly observe the linear convergence guaranteed by the theory in \Cref{section:linearconvergence}.
Although in \Cref{fig:lasso}, ADMM and NysADMM quickly exhibit linear convergence, \Cref{fig:elastic-net} clearly shows strong convexity leads to an improvement in the number of iterations required to converge.
Moreover, we see methods that make better use of curvature information converge faster than methods that do not, consistent with the results of the lasso and logistic regression experiments.  
\else
\fi

\begin{figure}[t] 
\captionsetup[sub]{font=scriptsize}
    \centering
    \begin{subfigure}[t]{0.46\textwidth}
        \centering
        \includegraphics[width=\columnwidth]{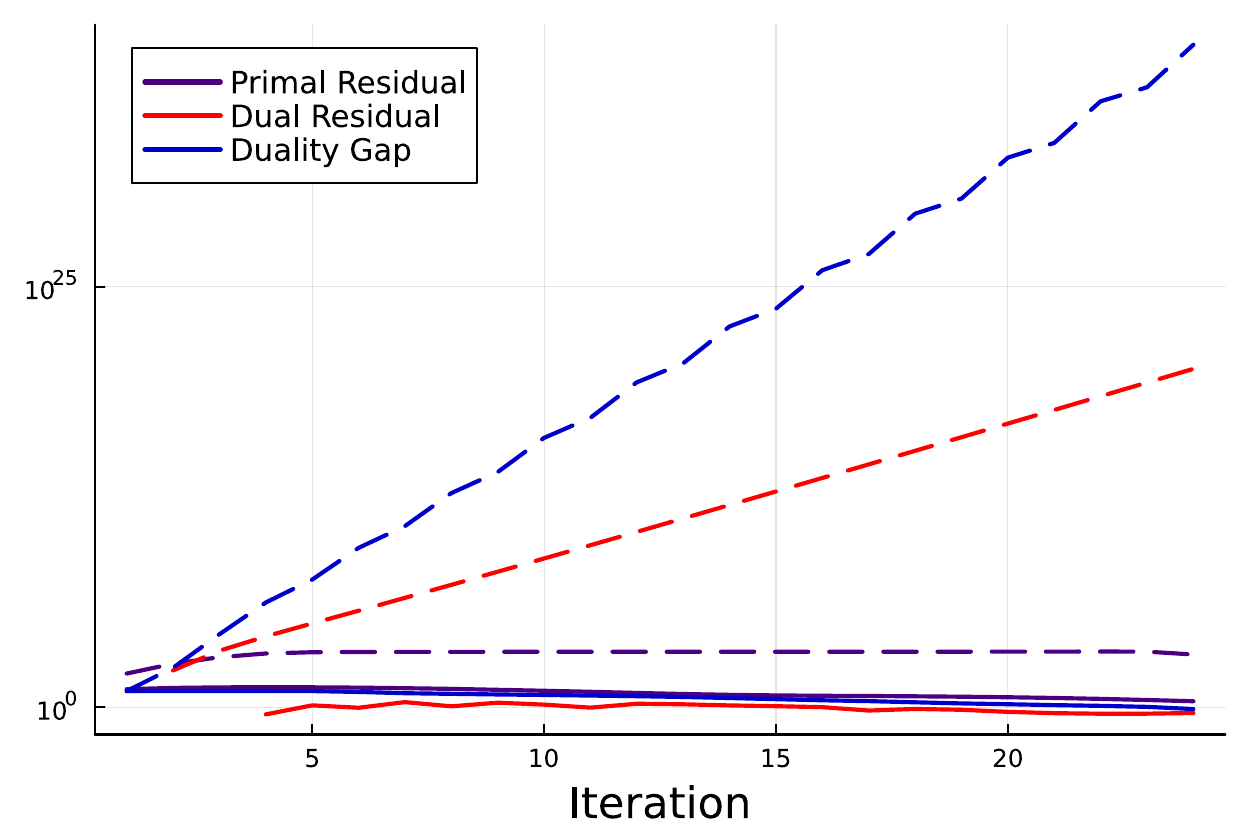}
        \caption{Lasso regression}
    \end{subfigure}
    \hfill
    \begin{subfigure}[t]{0.46\textwidth}
        \centering
        \includegraphics[width=\columnwidth]{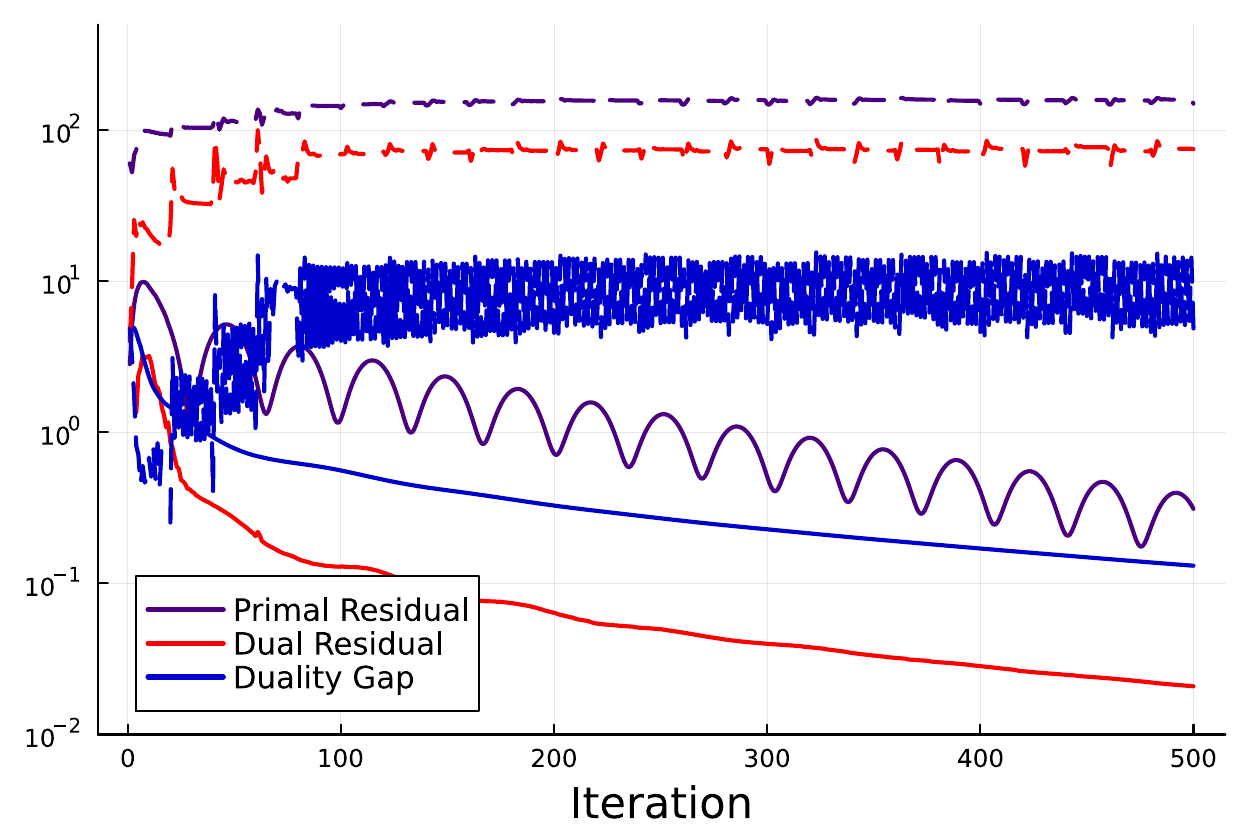}
        \caption{Logistic regression}
    \end{subfigure}
    \caption{We show the convergence or lack thereof for sketch-and-solve ADMM with (solid lines) and without (dashed lines) the correction term required by the theoretical results (see section~\ref{subsection:SketchAndSolveADMM}). When this term is not included, the algorithm does not converge.}
    \label{fig:divergence}
\end{figure}

\subsection{Sketch-and-solve ADMM fails to converge without the correction term}
\label{subsection:correction_exp}
To demonstrate the necessity of the correction term in \cref{subsection:SketchAndSolveADMM}, 
we run sketch-and-solve ADMM on lasso and $\ell_1$-logistic regression without the correction term.
\Cref{fig:divergence} presents the results of these simulations. 
Without the correction term in~\eqref{eq:sketch-and-solve}, sketch-and-solve ADMM quickly diverges on the lasso problem. For the logistic regression problem, it oscillates and fails to converge.
These results highlight the importance of selecting this term appropriately as discussed in \cref{subsection:SketchAndSolveADMM}.

\section{Conclusion and future work}
In this paper, we have developed a novel framework, GeNI-ADMM, that facilitates efficient theoretical analysis of approximate ADMM schemes and can aid the design of new, practical approximate ADMM methods. 
GeNI-ADMM generalizes prior approximate ADMM schemes as special cases by allowing various approximations to the $x$ and $z$-subproblems, which can be solved inexactly.
We have established the usual ergodic $\mathcal O(1/t)$-convergence rate for GeNI-ADMM under standard hypotheses, and linear convergence under strong convexity. 
We have shown how to derive explicit rates of convergence for ADMM variants that 
exploit randomized numerical linear algebra using the GeNI-ADMM framework.
Specifically, we have provided convergence results for NysADMM and sketch-and-solve ADMM, resolving whether these schemes are globally convergent. 
Numerical experiments on real-world data generally show an initial sublinear phase, followed by linear convergence, 
validating the theory we developed.

\section*{Acknowledgments}
The authors would like to thank Liubomir Baicev for carefully reading the manuscript and providing helpful comments. 
Z.\ Frangella, and M.\ Udell gratefully acknowledge support from
the National Science Foundation (NSF) Award IIS-2233762,
the Office of Naval Research (ONR) Award N000142212825 
and N000142312203, 
and the Alfred P. Sloan Foundation.
T.\ Diamandis is supported by the Department of Defense (DoD) through the National Defense Science \& Engineering Graduate (NDSEG) Fellowship Program.
B.\ Stellato is supported by the NSF CAREER Award ECCS-2239771.

\appendix
\section{Proofs not appearing in the main paper}

\ifpreprint
\subsection{Proofs for \cref{section:TechPrelimsAssump}}
We begin with the following lemma, which plays a key role in the proof of \cref{lemma:zsubprob}.
For a proof, see Example 1.2.2. in Chapter XI of \cite{hiriart1993convex}.
\begin{lemma}[$\varepsilon$-subdifferential of a quadratic function]
\label{lemma:QuadTecLemma}
Let $h(z) = \langle b,z\rangle +\frac{1}{2}\|z\|_A^2$, where $A$ is a symmetric positive definite linear operator. Then
    \begin{align*}
     \partial_{\varepsilon}h(z) &= \left\{b+A(z+v)\hspace{2pt}\middle \vert \hspace{2pt} \frac{\|v\|_A^2}{2} \leq \varepsilon\right\}.
     \end{align*}
\end{lemma}
With \cref{lemma:QuadTecLemma} in hand, we now prove \cref{lemma:zsubprob}.
\paragraph{Proof of \cref{lemma:zsubprob}}
\begin{proof} 
Observe the function defining the $z$-subproblem may be decomposed as $ G(z) =g(z)+h(z)$ with 
\[
h(z) = \frac{\rho}{2}\|M\tilde x^{k+1} +Nz + \tilde u^k\|^2_2+\frac{1}{2}\|z-\tilde z^k\|_{\Psi^k}^2.
\] 
Now, by hypothesis $\tilde{z}^{k+1}$ gives an $\varepsilon_z^k$-minimum of $G(z)$. 
Hence by \cref{proposition:epssubgrad-properties}, 
\[0\in \partial_{\varepsilon_z^k}G(\tilde z^{k+1})~\text{and}~\partial_{\varepsilon^k_z}G(\tilde z^{k+1}) \subset \partial_{\varepsilon^k_z}g(\tilde z^{k+1}) + \partial_{\varepsilon^k_z}h(\tilde z^{k+1}).\] 
Thus, we have $0 = s_g+s_h$ where $s\in \partial_{\varepsilon_z^k}g(\tilde z^{k+1})$ and $s_h \in \partial_{\varepsilon_z^k}h(\tilde z^{k+1})$.
Applying \cref{lemma:QuadTecLemma}, with $A = \Psi^k+\rho N^TN$ and $b = \rho N^{T}(M\tilde x^{k+1}+\tilde u^k)-\Psi^k\tilde z^k$,
we reach 
\[
s_h = \rho N^{T}(M\tilde x^{k+1}+N\tilde z^{k+1}+\tilde u^k)+\Psi^{k}(\tilde z^{k+1}-\tilde z^k)+v.
\]
The desired claim now immediately follows from using $s = -s_h$.
\end{proof}
\else
\fi

\subsection{Proofs for \cref{section:SubLinConvergence}}
\label{subsection:SubLinConvPfs}
\ifpreprint
\paragraph{Proof of \cref{prop:mu-incoherence}}
\begin{proof}
Using the eigendecomposition, we may decompose $H_f$ as  
\[
H_f = \lambda_{1}(H_f)\mathbf{v}\mathbf{v}^{T}+\sum^{d}_{i=2}\lambda_{i}(H_f)\mathbf{v}_i\mathbf{v}^{T}_i.
\]
So
\begin{align*}
  \|x^\star\|^2_{H_f}  &= \lambda_1(H_f)\langle \mathbf{v}, x^\star \rangle^2+\sum^{d}_{i=2}\lambda_{i}(H_f)\langle \mathbf{v}_i, x_\star \rangle^2 \\
  &\leq \lambda_{1}(H_f)\mu_\mathbf{v}\|x^\star\|^2+\lambda_{2}(H_f)\sum^{d}_{i=2}\langle \mathbf{v}_i,x^\star \rangle^2\\
  &\leq \left(\mu_\mathbf{v}+\frac{\lambda_{2}(H_f)}{\lambda_{1}(H_f)}\right)\lambda_{1}(H_f)\|x^\star\|^2.
\end{align*}
Recalling that $\sigma = \tau \lambda_1(H_f)$ with $\tau \in (0,1)$, we may put everything together to conclude
\begin{align*}
    \frac{(\lambda_{1}(H_f)-\sigma)\|x^\star\|^2}{\|x^\star\|^2_{H_f}}\geq (1-\tau) \frac{\lambda_{1}(H_f)\|x^\star\|^2}{\|x^\star\|^2_{H_f}} \geq (1-\tau) \left(\mu_\mathbf{v}+\frac{\lambda_{2}(H_f)}{\lambda_{1}(H_f)}\right)^{-1},
\end{align*}
proving the first claim. 
The second claim follows immediately from the first. 
\end{proof}
\else
\fi




\paragraph{Proof of \cref{prop:bnded_iterates}}
\begin{lemma}[Exact Gap function 1-step bound ]
Let $w^{k+1} = (x^{k+1}, z^{k+1}, u^{k+1})$ denote the iterates generated by \Cref{alg-gen_newton_admm} at iteration $k+1$, when we solve the subproblems exactly. 
Then, the following inequality holds
\begin{align*}
    	Q(w, w^{k+1}) & \le \frac{1}{2}\left(\|z-\tilde x^k\|^2_{\Theta^k}-\|x-x^{k+1}\|^2_{\Theta^k}\right)+\frac{1}{2}\left(\|z-\tilde z^k\|^2_{\Psi^k}-\|z-z^{k+1}\|^2_{\Psi^k}\right)\\
        &\frac{\rho}{2} \left(\|\tilde z^k - Mx\|^2 - \|z^{k+1} - Mx\|^2 \right)+\frac{\rho}{2}\left(\|\tilde u^k-u\|^2-\|u^{k+1}-u\|^2\right). 
\end{align*}
\end{lemma}
\proof
The argument is identical to the proof of \cref{lemma:GapBndInexactCase}.
Indeed, the only difference is $\varepsilon_x^k=\varepsilon_z^k =0$, as the subproblems are solved exactly.
So the same line of argumentation may be applied, except the exact subroblem optimality conditions are used. 
$\Box$

We now prove \cref{prop:bnded_iterates}.
\begin{proof}
First, observe that item 2. is an immediate consequence of item 1, so it suffices to show item 1.
To this end, plugging in $w = w^\star$, using $Q(w^\star,w^{k+1})\geq 0$, and rearranging, we reach
\begin{align*}
    &\|x^{k+1}-x^\star\|^2_{\Theta^k}+\|z^{k+1} - z^\star\|_{\Psi^k+\rho N^TN}^2+\rho\|u^{k+1}-u^\star\|^2\\
    &\leq \|\tilde x^{k}-x^\star\|^2_{\Theta^k}+\|\tilde z^k - z^\star\|_{\Psi^k+\rho N^TN}^2+\rho\|\tilde u^{k}-u^\star\|^2.
\end{align*}
Defining the norm $\|w\|_{W,k} = \sqrt{\|x\|^2_{\Theta^k}+\|z\|^2_{\Psi^k+\rho N^TN}+\rho \|u\|^2}$, the preceding inequality may be rewritten as 
\[
\|w^{k+1}-w^\star\|_{W,k}\leq \|\tilde w^k-w^\star\|_{W,k}.
\]
Now, our inexactness hypothesis (Assumption \ref{assp:x-subprob-orc}) along with $\nu$-strong convexity of the $z$-subproblem (which follows by Assumption \ref{assp:regularity}) implies
\[\|\tilde x^{k+1}-x^{k+1}\|\leq \varepsilon_x^k,\quad \|\tilde z^{k+1}-z^{k+1}\|\leq \sqrt{\frac{\varepsilon_z^k}{\nu}}.\]
Using $\tilde u^{k+1} = M\tilde x^{k+1}+N\tilde z^{k+1}+\tilde u^k$ and $u^{k+1} = Mx^{k+1}+Nz^{k+1}+\tilde u^k$, we find 
\[
\rho\|\tilde u^{k+1}-u^{k+1}\| \leq \varepsilon_x^k+C\sqrt{\varepsilon_z^k},
\]
where $C$ is some constant. 
Hence we have,
\[
\|\tilde w^{k+1}-w^\star\|_{W,k}\leq \|w^{k+1}-w^\star\|_{W,k}+\|\tilde w^{k+1}-w^{k+1}\|_{W,k}\leq  \|\tilde w^k-w^\star\|_{W,k}+C\left(\varepsilon_x^k+\sqrt{\varepsilon_z^k}\right).
\]
Defining the summable sequence $\varepsilon_w^k = C(\varepsilon_x^k+\sqrt{\varepsilon_z^k})$,
the preceding display may be written as 
\[
\|\tilde w^{k+1}-w^\star\|_{W,k}\leq \|\tilde w^k-w^\star\|_{W,k}+\varepsilon_w^k \quad (\triangle),
\]
Now, Assumption \ref{assp:metric_evol} implies that $\|\tilde w^k-w^\star\|_{W,k}\leq \sqrt{1+\zeta^k}\|\tilde w^k-w^\star\|_{W,k-1}$.
Combining this inequality with induction on $(\triangle)$, we find
\[
\|\tilde w^{k+1}-w^\star\|_{W,k}\leq \left(\prod^{k}_{j=2}\sqrt{1+\zeta^j}\right)\|\tilde w^1-w^\star\|_{W,1}+\left(\prod^{k}_{j=2}\sqrt{1+\zeta^j}\right)\sum_{j=1}^{k}\varepsilon_w^j.
\]
Hence,
\[
\sqrt{\min\{\theMin,\nu\}}\|\tilde w^{k+1}-w^\star\|\leq \|\tilde w^{k+1}-w^\star\|_{W,k}\leq \sqrt{\tau_\zeta} \left(\|\tilde w^1-w^\star\|_{W,1}+\mathcal E_w\right).
\]
It follows immediately that:
\[
\sup_{k\geq 1}\|\tilde w^{k+1}-w^\star\|\leq \max\left\{\frac{\sqrt{\tau_\zeta} \left(\|\tilde w^1-w^\star\|_{W,1}+\mathcal E_w\right)}{\sqrt{\min\{\theMin,\nu\}}},\|\tilde w^{0}-w^\star\|\right\}<\infty, 
\]
and so the sequence $\{\tilde w^{k}\}_{k\geq 0}$ is bounded, which in turn implies $\{\tilde x^{k}\}_{k\geq 0},\{\tilde z^{k}\}_{k\geq 0},$ and $\{\tilde u^{k}\}_{k\geq 0}$ are bounded.
\end{proof}

\paragraph{Proof of \cref{lemma:inexactxopt}}
\begin{proof}
Throughout the proof, we shall denote by $S_x^k(x)$, the function defining the $x$-subproblem at iteration $k$. 
The exact solution of the $x$-subproblem shall be denoted by $x^{k+1}$. 
To begin, observe that: 
\[
\nabla S^k_x(x) = \nabla f_1(x)+\nabla f_2(\tilde x^k)+(\Theta^k +\rho M^{T}M)(x-\tilde x^k)+\rho M^{T}(M\tilde x^k+N\tilde z^k+\tilde u^k).
\]
Now,
\begin{align*}
        \langle \nabla S^k_x(\tilde x^{k+1}),x^\star-\tilde x^{k+1}\rangle  & = \langle \nabla S^k_x(\tilde x^{k+1})-\nabla S^k_x(x^{k+1})+\nabla S^k_x(x^{k+1}),x^\star-x^{k+1}+x^{k+1}-\tilde x^{k+1}\rangle\\
        & = \langle \nabla S^k_x(\tilde x^{k+1})-\nabla S^k_x(x^{k+1}),x^\star-x^{k+1}\rangle\\
        &+\langle \nabla S^k_x(\tilde x^{k+1})-\nabla S^k_x(x^{k+1}),x^{k+1}-\tilde x^{k+1}\rangle\\
        &+\langle \nabla S^k_x(x^{k+1}),x^\star-x^{k+1}\rangle+\langle \nabla S^k_x(x^{k+1}),x^{k+1}-\tilde x^{k+1}\rangle\\
        &\overset{(1)}{\geq} \langle \nabla f_1(\tilde x^{k+1})-\nabla f_1(x^{k+1}),x^\star-\tilde x^{k+1}\rangle+\langle \tilde x^{k+1}-x^{k+1},x^\star-x^{k+1}\rangle_{\Theta^k+\rho M^{T}M}\\
        &+\langle \nabla f_1(\tilde x^{k+1})-\nabla f_1(x^{k+1}),\tilde x^{k+1}-x^{k+1}\rangle+ \langle \tilde x^{k+1}-x^{k+1},\tilde x^{k+1}-x^{k+1}\rangle_{\Theta^k+\rho M^{T}M}\\
        &+\langle \nabla S^k_x(x^{k+1}),x^{k+1}-\tilde x^{k+1}\rangle\\
        &\overset{(2)}{\geq} -C_1\varepsilon_x^k-C_2(\varepsilon^k_x)^2+\langle \nabla S^k_x(x^{k+1}),x^{k+1}-\tilde x^{k+1}\rangle.
\end{align*} 
Here $(1)$ uses $\nabla S_x^k(x)-\nabla S_x^k(y) = \nabla f_1(x)-\nabla f_1(y)+(\Theta^k+\rho M^{T}M)(x-y)$, and that $x^{k+1}$ is the exact solution of the $x$-subproblem,
while $(2)$ uses that $\tilde x^{k+1}$ is an $\varepsilon_x^k$ minimizer of the $x$-subproblem and that $f_1$ has a Lipschitz gradient. 
Now, as the iterates all belong to a compact set, it follows that $\sup_k \|\nabla S^k_x(x^{k+1})\| \leq C_3$.
So,
\begin{align*}
        \langle \nabla S^k_x(\tilde x^{k+1}),x^\star-\tilde x^{k+1}\rangle\geq -C_1\varepsilon_x^k-C_2(\varepsilon^k_x)^2 -C_3 \varepsilon_x^k \geq -C_x \varepsilon_x^k.
\end{align*}
This establishes the first inequality, and the second inequality follows from the first by plugging $\tilde x^{k+1}$ into the expression for $\nabla S_x^k(x)$.
\end{proof}

\paragraph{Proof of \cref{lemma:inexactzopt}}
\begin{proof}
By hypothesis, $\tilde z^{k+1}$ is an $\varepsilon_z^k$-approximate minimizer of the $z$-subproblem, so \cref{lemma:zsubprob} shows there exists a vector $\tilde s$ with $\|\tilde s\| \leq C\sqrt{\varepsilon_z^k}$, such that 
\[
g(z^\star)-g(\tilde z^{k+1})\geq \rho \langle \tilde u^{k+1}+\Psi^k(\tilde z^k-\tilde z^{k+1})+\tilde s,z^\star-\tilde z^{k+1}\rangle-\varepsilon_z^k.
\]
Rearranging, using Cauchy-Schwarz, along with the boundedness of the iterates, the preceding display becomes
\[
g(z^\star)-g(\tilde z^{k+1})- \rho \langle \tilde u^{k+1}+\Psi^k(\tilde z^k-\tilde z^{k+1}), z^\star-\tilde z^{k+1}\rangle\geq -C\sqrt{\varepsilon_z^k}-\varepsilon_z^k \geq -C_z\sqrt{\varepsilon_z^k},
\]
where the last display uses $\varepsilon_z^k \leq K_{\varepsilon_z}$, and absorbs constant terms to get the constant $C_z$. 
\end{proof}

\ifpreprint
\subsection{Proofs for \cref{section:linearconvergence}}
\label{section:linearconvergenceproof}
\paragraph{Proof of \cref{lemma:xzconvexity}}
\proof
\begin{enumerate}
\item We begin by observing that strong convexity of $f$ implies
    \begin{align*}
        \langle \nabla f(\tilde x^{k+1}) - \nabla f(x^\star), \tilde x^{k+1}-x^\star \rangle \geq \sigma_f\|\tilde x^{k+1}- x^\star\|^2.
    \end{align*}
    Decomposing $f$ as $f = f_1+f_2$, the preceding display may be rewritten as
    \begin{align*}
        &\langle \nabla f_1(\tilde x^{k+1})+\nabla f_2(\tilde x^{k}) - \nabla f(x^\star), \tilde x^{k+1}-x^\star \rangle + \langle \nabla f_2(\tilde x^{k+1}) - \nabla f_2(\tilde x^k), \tilde x^{k+1}-x^\star \rangle \\
        &\geq \sigma_f\|\tilde x^{k+1}- x^\star\|^2.
    \end{align*}
    Now, the exact solution $x^{k+1}$ of the $x$-subproblem at iteration $k$ satisfies $\nabla S_x^k(x^{k+1}) = 0$. 
    Moreover, $\tilde x^{k+1}$ is an $\varepsilon_x^k$ minimizer of the $x$-subproblem, and $\nabla S_x^k(x)$ is Lipschitz continuous.
    Combining these three properties, it follows there exists a vector $\Err$ satisfying: 
    \[
    \nabla S^k_x(\tilde x^{k+1}) = \nabla S_x^k(x^{k+1})+\Err = \Err, \quad \text{where}~\|\Err\|\leq C\varepsilon_x^k.
    \]
    Consequently, we have
    \[
    \nabla f_1(\tilde x^{k+1})+\nabla f_2(\tilde x^k) = \Theta^k(\tilde x^k-\tilde x^{k+1})+\rho M^{T}N(\tilde z^{k+1}-\tilde z^{k})-\rho M^{T}\tilde u^{k+1}+\Err.
    \]
    Utilizing the preceding relation, along with the fact that the iterates are bounded, we reach
    \begin{align*}
        & \langle \tilde x^k-\tilde x^{k+1},\tilde x^{k+1}-x^\star\rangle_{\Theta^k} + \rho\langle N(\tilde z^{k+1}-\tilde z^{k})+u^\star-\tilde u^{k+1},M(\tilde x^{k+1}-x^{\star})\rangle + C\varepsilon_x^k \\
        & \geq \sigma_f \|\tilde x^{k+1}-x^{\star}\|^2-\langle \nabla f_2(\tilde x^{k+1})-\nabla f_2(\tilde x^k),\tilde x^{k+1}-x^\star\rangle.~\Box
    \end{align*}
    \item As $-\rho N^{T}u^\star \in \partial g(z^\star)$, and $-\rho N^T\tilde u^{k+1}+\Psi^k(\tilde z^k-\tilde z^{k+1})+\tilde s \in \partial_{\varepsilon_z^k} g(\tilde z^{k+1})$, we have
    \begin{align*}
        &g(\tilde z^{k+1})-g(z^\star)\geq \langle -\rho N^T u^\star,\tilde z^{k+1}-z^\star\rangle, \\
        &g(z^\star)-g(\tilde z^{k+1})\geq \langle -\rho N^T\tilde u^{k+1}+\Psi^k(\tilde z^k-\tilde z^{k+1}),z^\star-\tilde z^{k+1}\rangle-\varepsilon_z^k.
    \end{align*}
    So, adding together the two inequalities and rearranging, we reach 
    \[
    \langle \tilde z^{k+1}-z^\star,\rho N^T(u^\star-\tilde u^{k+1})+\Psi^k(\tilde z^k-\tilde z^{k+1})+\tilde s\rangle \geq -\varepsilon_z^k.
    \]
    Now using boundedness of the iterates and $\|\tilde s\|\leq C\sqrt{\varepsilon_z^k}$, Cauchy-Schwarz yields
    \[
    \langle \tilde z^{k+1}-z^\star,\rho N^T(\tilde u^{k+1}-u^\star)+\Psi^k(\tilde z^k-\tilde z^{k+1})\rangle \geq -C\sqrt{\varepsilon_z^k}.
    \]

$\Box$
\end{enumerate}

\paragraph{Proof of \cref{lemma:expdescent}}
We wish to show for some $\delta>0$, that the following inequality holds
\[
\|w^k-w^\star\|_{G^k}^2-\|w^{k+1}-w^\star\|_{G^k}^2+\epsTot\geq \delta\|w^{k+1}-w^\star\|_{G^k}^2.
\]
To establish this result, it suffices to show the inequality
\[\|w^{k+1}-w^k\|^2_{G^k}+\frac{2\sigma_f}{\rho}\|x^{k+1}-x^\star\|^2-\frac{2}{\rho}\langle \nabla f_2(x^{k+1})-\nabla f_2(x^k),x^{k+1}-x^\star\rangle\geq \delta \|w^{k+1}-w_\star\|_{G^k}^2.
\]
We accomplish this by upper bounding each term that appears in $\|w^{k+1}-w_\star\|_{G^k}^2$, by terms that appear in the left-hand side of the preceding equality. 
To that end, we have the following result. 
\begin{lemma}[Coupling Lemma] \label{lemma:coupling_lemma} 
Under the assumptions of \cref{thm:linearconvabs}, the following statements hold:
    \begin{enumerate}
        \item Let $\mu_1\geq 2$, $c_1 = \frac{12L^2_f}{\rho^2\lamM},c_2 = \frac{4\left(2L^2_f+\theMax^2\right)}{\rho^2\lamM},c_3 = \frac{8\|M\|^2}{\lamM}$. Then for some constant $C\geq 0$, we have 
        \[
        \|u^{k+1}-u^\star\|^2\leq c_1\|\tilde x^{k+1}-x^\star\|^2+c_2\|\tilde x^{k}-\tilde x^{k+1}\|^2+c_3\|\tilde z^k-\tilde z^{k+1}\|_{1/\rho\Psi^k+N^TN}^2+C\epsTot.
        \]
        \item There exists constants $c_4 = 2\left(\frac{\|N\|^2+\psiMax/\rho}{\sigma^2_{\textrm{min}}(N)}\right)\|M\|^2$, $c_5 = c_4/\|M\|^2$, such that
        \[
        \|z^{k+1}-z^\star\|_{1/\rho \Psi^k+N^TN}^2\leq c_4\|x^{k+1}-x^\star\|^2+c_5\|u^{k+1}-u^k\|^2.
        \]
        \item For all $\mu>0$, we have 
        \begin{align*}
            \langle \nabla f_2(x^{k+1})-\nabla f_2(x^k), x^{k+1}-x^\star \rangle\leq \mu/2\|x^{k+1}-x^{\star}\|^2+\frac{L^2_f}{2\mu}\|x^{k+1}-x^{k}\|^2. 
        \end{align*}
    \end{enumerate}
\end{lemma}
Taking \cref{lemma:coupling_lemma} as given, let us prove \cref{lemma:expdescent}.
\proof
Observe \cref{lemma:coupling_lemma} implies 
\begin{align*}
    &\frac{1}{2}\|w^{k+1}-w^k\|^2_{G^k}+\frac{2\sigma_f}{\rho}\|x^{k+1}-x^\star\|^2-\frac{2}{\rho}\langle \nabla f_2(x^{k+1})-\nabla f_2(x^k),x^{k+1}-x^\star\rangle - \delta \|w^{k+1}-w_\star\|_{G_k}^2\\
    &\geq \left(\frac{2\sigma_f-\mu_2}{\rho}-\delta(c_1+c_4)\right)\|\tilde x^{k+1}-x^\star\|^2+\left(\frac{\theMin}{2\rho}-\frac{L_f^2}{2\rho\mu_2}-\delta c_2\right)\|\tilde x^{k}-\tilde x^{k+1}\|^2\\
    &+\left(1/2-\delta c_3\right)\|\tilde z^k-\tilde z^{k+1}\|_{1/\rho \Psi^k+N^TN}^2+\left(1/2-\delta c_5\right)\|\tilde u^k-\tilde u^{k+1}\|^2-C\epsTot.
\end{align*}
Setting $\mu_2 = \sigma_f$, and using $\theMin>L^2_f/\sigma_f = \kappa_f L_f$, we find by setting 
\[
\delta = \min\left\{\frac{\sigma_f}{\rho(c_1+c_4)}, \frac{\theMin-\kappa_f L_f}{2\rho c_2}, \frac{1}{2c_3},\frac{1}{2c_5}\right\}>0,
\]
that
\[
\|\tilde w^{k}-w^\star\|^2_{G^k}-\|\tilde w^{k+1}-w^\star\|^2_{G^k}+C\epsTot\geq \delta \|\tilde w^{k+1}-w^\star\|^2_{G^k},
\]
as desired. $\Box$

We now turn to to the proof of \cref{lemma:coupling_lemma}.
\paragraph{Proof of \cref{lemma:coupling_lemma}}
\proof
    \begin{enumerate}
        \item Observe by $L_f$-smoothness of $f$, that
        \begin{align*}
            2L^2_f\left(\|\tilde x^{k+1}-x^\star\|^2+\|\tilde x^k-x^{\star}\|^2\right) & \geq \|\nabla f_1(\tilde x^{k+1})-\nabla f_1(x^\star)+\nabla f_2(\tilde x^k)-\nabla f_2(x^\star)\|^2 \\
            & = \|\Theta^k(\tilde x^k-\tilde x^{k+1})+\rho M^{T}N(\tilde z^k-\tilde z^{k+1})+\rho M^{T}(u^{\star}-\tilde u^{k+1})+\Err\|^2\\
            &\geq \frac{1}{2}\rho^2\|M^{T}(\tilde u^{k+1}-u^\star)\|^2-\|\Theta^k(\tilde x^k-\tilde x^{k+1})+\rho M^{T}N(\tilde z^k-\tilde z^{k+1})+\Err\|^2.\\
        \end{align*}
       Here the last inequality uses with $\mu = 2$, the identity (valid for all $\mu>1$)
        \[\|a+b\|^2\geq (1-\mu^{-1})\|a\|^2+(1-\mu)\|b\|^2.\]
        So, rearranging and using $\|a+b\|^2 \leq 2\|a\|^2+2\|b\|^2$, we reach
        \begin{align*}
            \frac{\rho^2}{2}\|M^{T}(\tilde u^{k+1}-u^\star)\|^2 &\leq 6L^2_f\|\tilde x^{k+1}-x^{\star}\|^2+2\left(2L^2_f+\|\Theta^k\|^2\right)\|\tilde x^k-\tilde x^{k+1}\|^2\\
            &+ 4\rho^2 \|M\|^2\|N(\tilde z^k-\tilde z^{k+1})\|^2+4\|\Err\|^2.\\
        \end{align*}
     Now, using $\|\Err\|^2 \leq C\varepsilon_x^k\leq C\epsTot$ and $\|v\|_{N^TN}\leq \|v\|_{1/\rho \Psi^k+N^TN}$, we reach 
    \begin{align*}
        &\|\tilde u^{k+1}-u^\star\|^2\leq \frac{12 L^2_f}{\rho^2\lamM}\|\tilde x^{k+1}-x^{\star}\|^2+\frac{4\left(2L^2_f+\theMax^2\right)}{\rho^2\lamM}\|\tilde x^k-\tilde x^{k+1}\|^2\\
        & \frac{8\|M\|^2}{\lamM}\|\tilde z^k-\tilde z^{k+1}\|_{1/\rho \Psi^k+N^{T}N}^2+C\epsTot.~\Box
    \end{align*}
    \item This inequality is a straightforward consequence of the relation 
    \[
    M(\tilde x^{k+1}-x^\star) = \tilde u^{k+1}-\tilde u^k+ N(z^{\star}-\tilde z^{k+1}).
    \]
    Indeed, using the identity $\|a+b\|^2\leq 2\|a\|^2+2\|b\|^2$, we reach
    \begin{align*}
        & \|N(\tilde z^{k+1}-z^{\star})\|^2 \leq 2\|M\|^2\|\tilde x^{k+1}-x^\star\|+2\|\tilde u^k-\tilde u^{k+1}\|^2.
    \end{align*}
    Consequently
    \begin{align*}
        \|\tilde z^{k+1}-z^\star\|_{1/\rho \Psi^k+N^TN}^2 &\leq 2\left(\frac{\|N\|^2+\psiMax/\rho}{\sigma^2_{\textrm{min}}(N)}\right)\|M\|^2\|\tilde x^{k+1}-x^\star\|\\
        &+2\left(\frac{\|N\|^2+\psiMax/\rho}{\sigma^2_{\textrm{min}}(N)}\right)\|\tilde u^k-\tilde u^{k+1}\|^2,
    \end{align*}
    which is precisely the desired claim.~$\Box$

    \item Young's inequality implies for all $\mu>0$, that
    \[
    \langle \nabla f_2(\tilde x^{k+1})-\nabla f_2(\tilde x^k),\tilde x^{k+1}-x^\star\rangle \leq \frac{1}{2\mu}\|\nabla f_2(\tilde x^{k+1})-\nabla f_2(\tilde x^k)\|^2+\frac{\mu}{2}\|\tilde x^{k+1}-x^\star\|^2.
    \]
    The desired claim now follows from $L_f$-smoothness of $f_2$.~$\Box$ 
    \end{enumerate}
\else
\fi

\bibliographystyle{plainnat}
\bibliography{references}  






\end{document}